\documentclass[11pt]{article}

\usepackage{amssymb,amsmath,amsfonts,eurosym,geometry,ulem,graphicx,caption,color,setspace,sectsty,comment,footmisc,caption,pdflscape,subfigure,array,hyperref}
\usepackage[T1]{fontenc}
\usepackage[utf8]{inputenc}

\usepackage{lmodern} 

\usepackage{amsmath,amssymb,amsfonts}
\usepackage{natbib}
 \bibpunct[, ]{(}{)}{,}{a}{}{,}%
 \def\bibsep{\smallskipamount}%

\usepackage{rotating}
\usepackage{booktabs}
\usepackage{fancyvrb}
\usepackage{makecell}
\usepackage{algorithm}
\usepackage{algorithmic}
\usepackage{amsmath}
\usepackage{adjustbox}
\usepackage{subfigure}

\newtheorem{theorem}{Theorem}

\usepackage{multirow}
\usepackage{enumitem}
\usepackage{lipsum}
\newlength\mylen
\newlist{mycases}{enumerate}{1}
\setlist[mycases,1]{label=\textbf{Case~\arabic*.}, 
  labelwidth=\dimexpr-\mylen-\labelsep\relax,leftmargin=0pt,align=right}

\newenvironment{proof}[1][Proof]{\noindent\textbf{#1.} }{\ \rule{0.5em}{0.5em}}

\newcolumntype{L}[1]{>{\raggedright\let\newline\\arraybackslash\hspace{0pt}}m{#1}}
\newcolumntype{C}[1]{>{\centering\let\newline\\arraybackslash\hspace{0pt}}m{#1}}
\newcolumntype{R}[1]{>{\raggedleft\let\newline\\arraybackslash\hspace{0pt}}m{#1}}

\DeclareMathOperator*{\argmax}{argmax}
\DeclareMathOperator*{\argmin}{argmin}

\begin{document}

\begin{titlepage}
\title{Integrated Offline and Online Learning to Solve a Large Class of Scheduling Problems}
\author{Anbang Liu\thanks{Department of Automation, Tsinghua University, Beijing, China. Email: liuab19@tsinghua.org.cn} 
\and Zhi-Long Chen\thanks{Robert H. Smith School of Business, University of Maryland, College Park, MD, United States. Email: zlchen@umd.edu} 
\and Jinyang Jiang\thanks{Department of Management Science and Information Systems, Guanghua School of Management, Peking University, Beijing, China. Email: jinyang.jiang@stu.pku.edu.cn} 
\and Xi Chen \thanks{Department of Automation, Tsinghua University, Beijing, China. Email: bjchenxi@tsinghua.edu.cn}}
\date{\today}

\maketitle

\begin{abstract}
\noindent In this paper, we develop a unified machine learning (ML) approach to predict high-quality solutions for single-machine scheduling problems with a non-decreasing min-sum objective function with or without release times. Our ML approach is novel in three major aspects. First, our approach is developed  for the entire class of the aforementioned problems. To achieve this, we exploit the fact that the entire class of the problems considered can be formulated as a  time-indexed formulation in a unified manner. We develop a deep neural network (DNN) which uses  the cost parameters in the time-indexed formulation as the inputs to effectively predict a continuous solution to  this formulation, based on which a feasible discrete solution is easily constructed. 
The second novel aspect of our approach lies in how the DNN model is trained.
In view of the NP-hard nature of the problems, labels (i.e., optimal solutions) are hard to generate for training. To overcome this difficulty, we generate and utilize a set of special instances, for which optimal solutions can be found with little computational effort, to train the ML model offline.
The third novel idea we employ in our approach is that we develop an online single-instance learning approach to fine tune the parameters in the DNN for a given online instance, with the goal of generating an improved solution for the given instance. To this end, we develop a feasibility surrogate that approximates the objective value of a given instance as a continuous function of the outputs of the DNN, which then enables us to derive gradients and update the learnable parameters in the DNN. 
Numerical results show that our approach can efficiently generate high-quality solutions for a variety of single-machine scheduling min-sum problems with up to 1000 jobs.\\
\vspace{0in}\\
\noindent\textbf{Keywords:} Single-machine scheduling; Machine Learning; Supervised training; Online learning; Time-indexed formulation\\

\bigskip

\end{abstract}
\setcounter{page}{0}
\thispagestyle{empty}
\end{titlepage}

\pagebreak \newpage

\doublespacing

\section{Introduction}\label{Sec:Introduction}
Machine scheduling problems involve sequencing and scheduling a given set of jobs on a given set of machines to optimize an objective function subject to a given set of constraints. These problems frequently arise in manufacturing, service, and computer systems, and are among the most fundamental combinatorial optimization problems \citep{pinedo2012scheduling, blazewics-etal-2013}.  Single-machine scheduling problems are the most important and most extensively studied classes of machine scheduling problems. 
In such a problem, a set of jobs $\mathcal{N}=\{1, ..., n\}$ is processed on a single machine such that each job $j\in \mathcal{N}$ requires a processing time $p_j$, and the machine can process only one job at a time. The problem is to schedule the jobs, which is equivalent to assigning a starting time for each job, so as to optimize a certain objective function. Depending on the specific constraints and the objective function involved, a job $j$ may be associated with additional parameters such as an importance weight $w_j$, a release time $r_j$ (i.e., the time the job arrives at the system and becomes available), and a due date $d_j$ (i.e., the desired time by which the job is completed).

A wide range of objective functions has been studied in the literature \citep{pinedo2012scheduling}. In this paper, we focus on the objectives of the min-sum type, i.e., minimizing the sum of the costs incurred by individual jobs, which can be represented mathematically as 
$\sum_{j\in \mathcal{N}} z_j(C_j)$, where $C_j$ is the completion time of job $j$ and $z_j$ is the cost function for job $j$ which is usually non-decreasing. Here the structure of the overall cost function $\sum_{j\in \mathcal{N}} z_j(C_j)$ is called {\it additive} or {\it sum-type} because it is the summation of the individual jobs' costs. 
Sum-type objective functions are among the most common categories of scheduling criteria studied in the scheduling literature. They include total weighted completion time $\sum_{j\in \mathcal{N}} w_jC_j$ and total weighted tardiness $\sum_{j\in \mathcal{N}} w_jT_j$, where $T_j = \max(0, C_j-d_j)$ defines the tardiness of job $j$, among others. A class of multi-criterion scheduling problems \citep{t2001multicriteria, minella2008review} involve minimization of a weighted sum of two or more additive functions,  e.g., 
$\rho \sum_{j\in \mathcal{N}}w_{1j}T_j  + (1-\rho)\sum_{j\in \mathcal{N}}w_{2j}C_j $, where $0<\rho <1$ is a given constant, and $w_{1j}$ and $w_{2j}$ are the importance weights of job $j$ for the two criteria, respectively. Such objective 
functions can be rewritten as sum-type functions, e.g., $\sum_{j\in \mathcal{N}} [(\rho w_{1j})T_j + ((1-\rho)w_{2j})C_j]$. Thus, those multi-criterion scheduling problems also have min-sum objectives.  In all the above-discussed min-sum objectives, the cost functions for different jobs, i.e., $z_j(C_j)$'s, have the same structure.
However, there are more complex min-sum objectives where different jobs may involve different cost function forms. A class of multi-agent scheduling problems \citep{agnetis2014multiagent} has such objectives. A multi-agent scheduling problem involves multiple agents, where each agent owns a subset of jobs and has a specific objective function to optimize for its jobs, while all the jobs are processed on the same machine. 
In one type of multi-agent scheduling problems, the objective is to minimize a weighted sum of a sum-type cost function of each agent, e.g.,  $\rho \sum_{j\in \mathcal{N}_A} w_{j}T_j + (1-\rho)\sum_{j\in \mathcal{N}_B} w_{j}C_j$, where there are two agents, each owning a subset of jobs $\mathcal{N}_A$ and $\mathcal{N}_B$, respectively. Such cost functions are obviously also of a sum type and can be written as  $\sum_{j\in \mathcal{N}} z_j(C_j)$, where $\mathcal{N} = \mathcal{N}_A\cup \mathcal{N}_B$, but different jobs may have different cost function forms, e.g., $z_j(C_j) = \rho w_{j}T_j$ for $j\in \mathcal{N}_A$ and $z_j(C_j) = (1-\rho)w_{j}C_j$ for $j\in \mathcal{N}_B$.  
Besides the above discussed objective functions, which are all piece-wise linear, we also consider more general nonlinear cost functions, including exponential function $z_j(C_j) = w_{j}C_j^{a_{j}}$  \citep{szwarc1988single,janiak2009scheduling}, where the exponent $a_{j}>0$ is job dependent. 
The approach we develop works for any type of min-sum objective as long as the sum-type objective function is non-decreasing in the completion times of the jobs. 

We adopt the commonly used three-field notation $\alpha |\beta | \gamma$ proposed by \cite{graham1979optimization} to represent a scheduling problem, where the $\alpha$ field describes the machine environment, the $\beta$ field describes the constraints, and the $\gamma$ field is the objective function to be minimized. In this paper, we consider the entire class of single-machine scheduling problems with a min-sum objective that is non-decreasing in job completion times, including those  without release times, i.e., $1||\sum z_j(C_j)$, where all the jobs are ready at time 0, and those with release times, i.e., $1|r_j|\sum z_j(C_j)$, where different jobs have generally different ready times. Most of the specific min-sum problems within the class of problems we study are known to be strongly NP-hard and hence are among the most difficult classes of combinatorial optimization problems \citep{garey-johnson-1979}. Problems with release times $1|r_j|\sum z_j(C_j)$ with any commonly studied cost functions $z_j$ are all strongly NP-hard because the simplest among them, which is problem $1|r_j|\sum C_j$, is already strongly NP-hard \citep{lenstra1977complexity}. Most single-machine problems without release times, $1||\sum z_j(C_j)$, are also strongly NP-hard except a handful of them that are solvable either in polynomial-time or in pseudo-polynomial time, including $1||\sum w_jC_j$, $1||\sum U_j$, $1||\sum w_jU_j$, and $1||\sum T_j$, where $U_j=1$ if $C_j> d_j$ and 0 otherwise. Problem  $1||\sum w_jC_j$, can be solved by sequencing the jobs in non-increasing ratios $w_j/p_j$ \citep{smith-1956}; problem $1||\sum U_j$ can be solved in polynomial time by Moore's algorithm \citep{moore-1968}; problems $1||\sum w_jU_j$ and $1||\sum T_j$ are ordinarily NP-hard but can be solved in pseudo-polynomial time \citep{sahni-1976,lawler-1977}. 

A variety of local search and mathematical programming-based heuristic algorithms exist for various NP-hard machine scheduling problems \citep{anderson1997machine,grimes2015solving, ɖurasevic2023heuristic}. 
To make such an approach effective, one often needs to exploit special structures and solution properties associated with the problem and customize the solution approach for the problem. However, inherent solution structures associated with one objective function (e.g., $\sum w_jC_j$) can be very different from those associated with another objective function (e.g., $\sum w_jT_j$). Consequently, approaches developed for one problem may not work well for another problem. Another drawback of such approaches is that they are often time-consuming, especially for instances with large sizes.  

The goal of this paper is to develop a unified machine learning approach that is effective for all the problems $1||\sum z_j(C_j)$ and $1|r_j|\sum z_j(C_j)$, and once trained, takes little time to solve an instance. Our idea is to formulate all these problems as a unified time-indexed binary integer program (BIP), and  develop a deep learning method based on this formulation. To make  it less time consuming to train the model on large instances, we develop a technique to  construct special large instances for which  optimal solutions can be easily found. 
Moreover, we integrate supervised offline training  and online single-instance learning to improve the solution quality for a given instance. We describe below in subsection~\ref{sec:TI-formulation} how our problems can be formulated as a time-indexed formulation, and give a brief introduction about the machine learning approaches we use in subsection~\ref{sec:machine-learning}. In subsection~\ref{sec:Contribution}, we summarize the major challenges and our contributions. 

\subsection{Time-indexed formulation}\label{sec:TI-formulation}

Time-indexed formulations for machine scheduling problems are widely adopted for formulating and solving machine scheduling problems \citep{sousa1992time,van2000time}. To create such a formulation, one needs to discretize the planning horizon, define time-indexed binary decision variables to represent possible starting times of the jobs, and define the objective function and constraints accordingly. All the problems considered in this paper can be formulated into a unified time-indexed BIP formulation as follows. Define the following parameters: 
\begin{itemize}
    \item Processing times of the jobs $\{p_j\}_{j\in \mathcal{N}}$, and the total processing time of the jobs, $P=\sum_{j\in \mathcal{N}} p_j$.
    \item Release times of the jobs $\{r_j\}_{j\in \mathcal{N}}$, which are zero for a problem without release times.
    \item Let $T$ be the length of the planning horizon, where the value of $T$ should be large enough such that in a schedule  without inserted idle time, all the jobs can be completed  by time $T$. Discretize the planning horizon into $T$ time points, $\mathcal{T}=\{0,1,...,T-1\}$. 
    \item Cost of job $j\in \mathcal{N}$ if it starts at time $t\in \mathcal{T}$, $c_{jt} = z_j(t+p_j)$. For ease of presentation, we call  $c_{jt}$ the {\it starting cost} of job $j$ if it is started at time $t$. Then, $\{c_{j0}, c_{j1}, \cdots, c_{j,T-1}\}$ together are called the {\it starting costs} of job $j$.
\end{itemize}  
For each job $j\in\mathcal{N}$ and each time point $t\in\mathcal{T}$, define a binary decision variable $x_{jt}$ to be 1 if job $j$ starts at time  $t$, and zero otherwise. Then we have the following unified time-indexed BIP formulation for all the problems we consider: 
\begin{subequations}
\label{time-index-formulation}
\begin{align}
    {\text{minimize}} &\
    \sum_{j\in \mathcal{N}}\sum_{t=0}^{T-p_j}c_{jt}\cdot x_{jt} \label{obj:timeindex}\\
    \text{subject to} &\
    \sum_{t=r_j}^{T-p_j}x_{jt}=1, ~~\forall ~j \in \mathcal{N}, \label{const:assignment}\\
    &\ \sum_{j\in \mathcal{N}}\sum_{k=\max(r_j, t-p_j+1)}^{t}x_{jk}\leq 1, ~~\forall ~t \in\mathcal{T}, \label{const:capacity}\\
    &\ x_{jt}\in \{0, 1\}, ~~\forall ~j \in \mathcal{N}, t\in\mathcal{T}. 
    \end{align}
\end{subequations}
Constraints \eqref{const:assignment} guarantee that each job has a unique starting time no earlier than its release time. The capacity constraints \eqref{const:capacity} ensure that at most one job is processed at a time.

A major advantage of this time-indexed formulation is that this single formulation represents all the single-machine scheduling problems that we intend to solve. Therefore, we build on this formulation to develop a unified machine learning approach for all these problems. Another advantage of such a time-indexed formulation is that its LP relaxation tends to generate tighter bounds than other commonly used integer programming formulations. However, such a formulation contains a large number of binary decision variables and constraints, especially when the total processing time of the jobs $P$ is large. Consequently, to solve such formulations, valid inequalities and other solution techniques such as column generation approaches are often used \citep{berghman-2015, van2000time}.

\subsection{Machine learning}\label{sec:machine-learning}
Many studies have developed machine learning (ML) methods to solve difficult combinatorial optimization problems, including scheduling problems and routing problems \citep{bengio2021machine}. In some studies \citep{morabit2023machine,zhang2022learning}, ML is integrated within optimization algorithms  to accelerate the computation. In this paper, we develop a supervised  learning method to directly generate a solution for a given instance. In the following, we  summarize the key ideas of supervised learning. 

Supervised learning is a major sub-field of machine learning \citep{lecun2015deep}. A supervised learning model can be viewed as a parameterized function $f_\theta$ with the learnable parameters $\theta$ that maps an input $x\in \mathcal{X}$ to an output $y\in \mathcal{Y}$, where $\mathcal{X}$ is the input space consisting of all possible values of the input, and $\mathcal{Y}$ is the output space consisting of all possible values of output. A supervised learning model $f_\theta$ learns  from existing data and generalizes to unseen data. Thus, a training data set, $\mathcal{D}^{\text{train}}=\{(x_1,y_1), (x_2, y_2), ... \}$, is required, where $x_i\in \mathcal{X}$ is an input and $y_i\in \mathcal{Y}$ is the corresponding correct output (also called label). The learning process is to find the optimal learnable parameters $\theta$  to minimize the average loss as follows,
\begin{equation}
\label{intro_equation_SL}
    \min_{\theta}  \frac{1}{|\mathcal{D}^{\mathrm{train}}|} \sum_{(x,y) \in \mathcal{D}^{\mathrm{train}}} \mathcal{L}(y,f_\theta(x)).
\end{equation}
where $\mathcal{L}(\cdot)$ is a problem-specific loss function. The trained model is used to generate outputs for unseen inputs. A model is said to have good generalization capability if it produces high-quality outputs for unseen inputs. Usually, the values in the input space  follow some unknown distribution  $\mathbb{P}$. A key to achieving good generalization is that the inputs in the training set be sampled independently and identically from the same distribution $\mathbb{P}$. Moreover, the training set should be sufficiently large to capture the properties of the input space.

 In the context of scheduling problems, an input $x$ is designed to be able to fully describe a given problem instance, and the input space $\mathcal{X}$ should cover all  possible instances to be solved online. The output $y$ is a solution of the instance. For example, for an instance of $1||\sum w_j T_j$, the input $x$ can consist of  job processing times, weights, and due dates, and $y$ can be the optimal job starting times. 
Training of a supervised learning model for difficult  scheduling problems can be time-consuming, but after the training is completed, the trained model  is often computationally efficient, with the time complexity being polynomial in the size of the instance.

There are various ways to design a learning model $f_{\theta}$. It can be designed as a relatively simple model, e.g., a linear model in Support Vector Machine approaches \citep{hearst1998support} or a combination of a linear model and a Sigmoid function in logistic regression \citep{lavalley2008logistic}. These models are typically easy to train. However, limited by the low learning ability, these models are often not capable of discovering complex features on their own. To have a good performance, features need to be carefully defined and selected. Deep Neural Networks (DNNs) are advanced machine learning models that have recently achieved success in image and speech processing. The key to the success lies in their depth. 
Unlike simpler models, a DNN model often consists of  a large number of layers of interconnected neurons and nonlinear activation functions.  The input can be raw data without feature engineering. Through training, the DNN learns feature representations automatically, with earlier layers capturing simpler features and deeper layers extracting more complex features.

\subsection{Challenges and contributions}
\label{sec:Contribution}
To develop a unified machine learning approach that works for all the  single-machine scheduling problems with a non-decreasing min-sum objective, there are two major challenges.

\begin{itemize}
    \item Different scheduling problems can exhibit very  different characteristics due to the wide range of sum-type objective functions. For example, $\sum w_jC_j$ is a smooth function of the completion times $C_j$, whereas $\sum w_jT_j$ is not. Moreover, in problems without release times, jobs are all available at time 0, whereas in  problems with release times, jobs may not be available before certain time points. It is thus a challenge to develop a unified machine learning approach that works well for  a wide variety of problems.
    \item As discussed above, to train a supervised learning model with a good generalization capability, a training set with sufficiently many instances  drawn i.i.d. from the distribution of the input space is needed.  However, because the scheduling problems under consideration are NP-hard, it is impractical to generate optimal solutions as labels for large-sized problems. Thus, once a machine learning model is developed for our problems, it is a challenge to train the model. As shown in the literature, if small-sized instances, which can be optimally solved, are used in the training set, the trained model usually does not have good generalization performance. Although reinforcement learning can be employed without the need for labels, its performance often deteriorates on large problems due to the extensive action space involved. The related literature is reviewed  in Section \ref{Sec:Literature_Review}.
\end{itemize}
In this paper, we address these two difficulties in three novel ways:

\begin{itemize}
\item  We design a unified machine scheduling neural network (UMSNN) by exploiting the time-indexed formulation and by leveraging the capability of the underlying deep neural networks (DNN) to automatically extract hidden features that determine optimal solutions. Specifically, we define the inputs of the UMSNN by using the time-indexed formulation. We use the starting costs of jobs, i.e., $\{c_{jt}\}_{j\in \mathcal{N}, t\in \mathcal{T}}$ as part of the inputs of our machine learning model. Since the objective function information is incorporated into the values of the starting costs, using the starting costs as the inputs makes the approach viable for all  min-sum problems. The other part of the inputs of the ML model includes job processing times and the job release times. Since $n$ and $T$ can be very large, the inputs defined in this manner are of high dimension, and are therefore unstructured and ``raw". To extract the hidden features that determine optimal solutions, the UMSNN is designed to have a deep structure. 
A cross-entropy loss function, which measures the distance between predictions and labels, is used.  

\item  We develop a novel offline training scheme using specially designed instances. Unlike existing methods that use  training instances with a smaller number of jobs  than in the testing instances, we train on large-sized instances which have a similar number of jobs as in the testing instances. However, the large training instances we use are specially constructed from  instances with much shorter processing times than those in the testing instances such that optimal solutions to such  instances can be found quickly and these solutions can be easily adapted to become the optimal solutions for the corresponding large instances with desired job processing times.

\item 
To improve the performance of the trained UMSNN for any given instance, we further develop  an online  learning approach for the given instance. Using the given instance as the input, and the objective function of the instance as the loss function to be minimized, the learnable parameters of the UMSNN are further optimized to minimize the loss function for this specific instance. Since the objective value as a function of  the UMSNN's outputs is non-differentiable, a feasibility surrogate is developed to approximate the objective value as a differentiable function of the outputs from UMSNN. The surrogate  provides reasonable gradients to allow for loss back-propagation.

\end{itemize}
 
To demonstrate the performance of our approach, we test both the approach with offline learning only, and the integrated approach with both offline and online single-instance learning, based on  randomly generated instances with 500 to  1000 jobs. It shows that while both approaches can generate high-quality solutions quickly,  the integrated approach generates even better solutions with slightly more computational effort than the approach with offline learning alone. For comparison purposes, we also test two benchmark approaches which also rely on the time-indexed formulation but do not use machine learning.

This paper is organized as follows. In Section \ref{Sec:Literature_Review}, existing machine learning approaches for scheduling  and related problems are reviewed. Section \ref{sec:UMSNN} develops the architecture of UMSNN. Section~\ref{sec:generation_of_trainingset} describes the offline supervised training using specially generated instances. In Section \ref{sec:online_learning}, the online single-instance learning approach is developed. Computational results are presented in Section \ref{Sec:Numerical_Results}, followed by concluding remarks in Section \ref{sec:conclusion}.

\section{Literature Review} \label{Sec:Literature_Review}
In this section, we first review supervised learning approaches for solving machine scheduling and related problems. These approaches require labels for training. Then, we review relevant studies using reinforcement learning approaches, where labels are not needed for training. Finally, we review other heuristic approaches that do not use  machine learning for some of the scheduling problems that we study.

\subsection{Supervised learning approaches}

A supervised learning model is trained by using labeled data, where the input data is paired with the correct output. The model learns to make predictions or decisions based on this training, aiming to generalize and accurately predict outcomes for new, unseen data. 
The existing approaches for scheduling problems can be categorized into two types, ``direct approaches" and ``indirect approaches". 
In direct approaches, which are also referred to as ``end-to-end" approaches, machine learning models directly learn to generate solutions for given problem instances. 
By contrast, in indirect approaches, machine learning methods are combined with mathematical methods or heuristic approaches to replace time-consuming components or to provide guidance for the search process. In the following, direct approaches are reviewed first, followed by indirect approaches. 

\textbf{End-to-end approaches:}  \cite{parmentier2023structured} study the single-machine scheduling problem $1|r_j|\sum C_j$, a strongly NP-hard problem, by developing a structured learning model. A total of 27 features (such as $\frac{p_j}{\sum_j p_j}$) are manually defined to capture the information for a given instance of $1|r_j|\sum C_j$ and serve as the inputs of the ML model. The outputs are the job processing times of an instance of $1||\sum C_j$, which is polynomially solvable with the optimal schedule being the non-decreasing order of job processing times. Therefore, essentially, the learning model learns to directly output the optimal schedule. 
The model is trained in a supervised manner. Due to the NP-hard nature of $1|r_j|\sum C_j$, it is impractical to generate optimal schedules as labels when problem sizes are large. 
\cite{parmentier2023structured} solve small-sized instances (where there are 50 to 110 jobs) optimally and use them as labels for supervised training. The trained model is generalized for solving unseen large-scale instances.  
As reported in \cite{parmentier2023structured}, the approach has the advantage in terms of computational efficiency.
However, the features are manually defined based on the properties of $1|r_j|\sum C_j$, limiting their applicability. For problems with a different objective function, a different set of features may have to be defined.

\cite{weckman2008neural} consider the job shop scheduling makespan problem by developing a multi-layer perceptron which is trained by using labels generated by genetic algorithms. The approach is tested on a set of small instances (with no more than 20 jobs and 20 machines). The results obtained are not obviously better than those obtained by a genetic algorithm. 
\cite{schmidt2021approaching} propose a hybrid supervised learning approach combining deep neural network and greedy heuristics for a parallel machine scheduling problem with sequence-dependent setup times. Makespan and total lateness are considered as objective functions. The commercial CPLEX solver is used to generate labels for supervised training. As reported in the paper, the performance of their approach is not obviously better than those obtained by priority dispatching rules.

In addition to scheduling problems, the traveling salesman problem (TSP) can also be solved by using end-to-end supervised learning approaches. In \cite{vinyals2015pointer}, a pointer network is developed based on the encoder-decoder Recurrent Neural Network to solve TSP. For a given instance, the input features are defined as the sequence of the city coordinates, and the output is a permutation of the cities, which indicates the order of visiting the cities. 
For small examples with no more than 20 cities, a set of instances and the corresponding optimal solutions are used for supervised training. High-quality solutions are predicted by the pointer network. 
However, for examples containing 20 to 50 cities, it is hard to generate optimal solutions as labels because of the NP-hardness of the problem.
Instead, they use solutions generated by approximation algorithms as labels for training. 
Since the labels do not have high quality, the solution qualities obtained are not better than those obtained by the approximation algorithms. 

\textbf{Indirect approaches:} 
Many heuristics and optimization approaches solve scheduling and related problems by searching the feasible solution space. When the solution space is large, such approaches are generally time-consuming. To reduce the computational efforts required, machine learning has been developed to learn to guide the search in the solution space. In \citet{bouvska2023deep}, single-machine scheduling problem  $1||\sum T_j$ is considered. The problem is ordinarily NP-hard, and can be solved by using Lawler’s decomposition \citep{lawler1977pseudopolynomial} and a symmetric decomposition proposed by \citet{della1998new} in Pseudo-polynomial time. \citet{bouvska2023deep} develop a DNN and integrate it within the two decomposition approaches to guide the decomposition by estimating the objective function related to a decomposition action. 
The DNN is trained in a supervised manner. To generate a training set, a set of problem instances is randomly generated, and is solved by using the two decomposition approaches without a DNN. Since each instance is recursively decomposed, a large number of decomposition actions and the corresponding objectives (labels) can be gathered while solving each instance. As reported in the paper, the proposed approach outperforms the state-of-the-art approaches \citep{garraffa2018exact, holsenback1992heuristic}.
The limitation of the approach is that it relies on the specific structure of $1||\sum T_j$. Thus, for  problems with a different objective function, this approach does not work. 

Besides using ML to guide the search of the solution space, there are also approaches that replace a time-consuming part of an optimization method by a supervised learning approach to save computational time.
In \cite{liu2023integrating}, job shop scheduling with the objective  of minimizing the total weighted tardiness is considered, and ML is integrated within Lagrangian Relaxation (LR) to solve subproblems. Their approach decomposes the overall scheduling problem into a set of subproblems by relaxing machine capacity constraints, and develops a DNN approach to solve subproblems in a computationally efficient manner. 
Since subproblems can be optimally solved in a reasonable time, a large number of subproblem instances and the corresponding labels can be generated in a reasonable time for supervised training. Their approach can find a solution to the overall scheduling problem  efficiently. However, the quality of the overall solution obtained depends on the tightness of the Lagrangian relaxation bound.

In addition to scheduling problems, many other optimization problems are solved in the literature by approaches that combine optimization and supervised learning.  
For example, in \cite{khalil2016learning,balcan2018learning,paulus2024learning}, ML is integrated within the branch-and-bound solution methodology for  optimization problems with integer decision variables, where ML is trained to select variables to branch on. 
ML can also be integrated within column generation approaches to select columns to add \citep{morabit2021machine,morabit2023machine}, or  learn to solve subproblems \citep{shen2022enhancing, vaclavik2018accelerating, minaeva2016scalable, burke2014new}.

\subsection{Reinforcement learning approaches}

One of the major difficulties associated with supervised learning approaches for solving NP-hard combinatorial optimization problems is that labels are difficult to generate. By contrast, the training process of a reinforcement learning (RL) approach is a process of trial and error with no labels needed. 
Using RL, problems to be solved are first formulated as Markov decision process (MDP) problems with multiple stages of decision making. At any stage $t$, the complete situation is encapsulated in a state $s_t$, and an action is taken. The state $s_t$ transitions to $s_{t+1}$ after taking the action, and a reward is received, which is problem-dependent. An RL approach is to train a model to generate an action for any state such that the total expected reward is maximized. Based on the way the actions are generated, RL approaches can be further categorized into value-based approaches, e.g., Q-learning, and policy-based approaches, e.g., policy gradient approaches.

There are RL-based approaches developed in the literature to solve scheduling problems \citep{wang2005application,yuan2016dynamic,zhang2020learning,li2023two}. In \cite{wang2005application}, an RL approach is developed for single-machine scheduling with the objective of minimizing the maximum lateness, the total number of tardy jobs, or the total tardiness. The approach incorporates multiple dispatching rules, and the action at a stage is to select a dispatching rule. A Q-learning algorithm is developed to learn to generate an action at each stage. In \cite{yuan2016dynamic}, dynamic parallel machine scheduling with breakdowns is considered and a Q-learning approach is developed to learn to select dispatching rules. The advantage of these methods is that they maintain the high computational efficiency of dispatching rules and allow for real-time scheduling. However, the quality of the generated solutions is questionable.

In several other studies, actions are defined to directly generate schedules of jobs. \cite{li2023two} consider parallel machine scheduling problems with family setups and the objective of minimizing the total tardiness. Whenever the machine is idle, an action is taken, which involves selecting a job from the available jobs and assigning it to the idle machine.
A proximal policy optimization (PPO) algorithm is developed to learn a good policy. 
In \cite{zhang2020learning}, job shop scheduling problems with the makespan objective are considered. Each stage is associated with a time slot, and the action at a stage is to start a job or not. A policy gradient approach is developed to learn a policy, which can generate an action based on a given state. 
Once trained, these approaches can handle large-scale problems within a reasonable time. The solution qualities are generally better than those obtained by dispatching rules. However, compared to the solutions  generated by optimization based approaches, the solutions generated by RL-based approaches are often inferior. As reported in \cite{zhang2020learning}, for many instances with less than 100 jobs, the gaps between the objective values obtained by their RL-based approach and the best-known values are over 30\%. The reason for the poor solution quality is that when the problem size is large, the feasible solution space of the problem becomes vast, and the action space of the MDP problem also becomes extensive, making it difficult for RL-based algorithms to learn high-quality strategies. Furthermore, the reward function needs to be defined based on the specific objective function considered. To the best of our knowledge, there is no unified RL approach that can handle various objective functions. 

Besides solving scheduling problems, in the literature, RL is also developed for solving many other problems such as traveling salesman problems \citep{kool2018attention}, vehicle routing problems \citep{nazari2018reinforcement}, and object packing  problems \citep{huang2022planning}.

\subsection{Heuristic approaches}

In this subsection, heuristic approaches for single machine scheduling problems with min-sum objectives are reviewed. The problem $1||\sum w_jT_j$ has been extensively studied in the literature.
In \citep{tasgetiren2004particle}, a particle swarm optimization algorithm is developed, successfully solving instances with up to 100 jobs in under 100 seconds. \citet{bilge2007tabu} proposes a Tabu search algorithm that obtains optimal solutions for instances with 40 jobs within 100 seconds. For instances with 100 jobs, the algorithm achieves solution qualities that surpass the best-known results, although the computational time increases to between 200 and 300 seconds. \citet{ding2016breakout} proposes a breakout dynasearch algorithm, solving instances with up to 300 jobs. For all considered instances, the algorithm is run 20 times, and a solution matching the optimal objective value is consistently obtained within an average of 252 seconds across these runs.
To enhance computational efficiency for larger instances, \citet{uchronski2021parallel} develops a parallel Tabu search algorithm, solving instances with up to 1,000 jobs. Compared to sequential algorithms, the proposed parallel approach offers significant improvements in computational efficiency.  However, the solution quality is questionable. The study uses the heuristic method proposed by \citep{potts1991single} as a benchmark, and for some instances, the objective values obtained exhibit a gap exceeding 10\% compared to the benchmark.
In \citep{cheung2017primal}, a linear programming (LP) relaxation-based pseudo-polynomial-time heuristic is developed for the entire class of min-sum single-machine scheduling problems without release dates. It is demonstrated that the proposed heuristic can generate a schedule with a cost at most four times the optimal value.

Single-machine scheduling problems with min-sum objectives and release dates are known to be strongly NP-hard. \citet{goemans2002single} examines the LP relaxations for the problem $1|r_j|\sum w_jC_j$, showing that the LP relaxations of both the time-indexed and completion-time formulations can produce solutions with errors of at most 1.6853.
\citet{avella2005near} develops a Lagrangian relaxation-based approach for $1|r_j|\sum w_jC_j$. The largest set of instances solved in this study includes 400 jobs, each with a maximum processing time of 50. Using approximately 2,200 seconds, the approach achieves optimality gaps ranging from 0.4\% to 0.9\%.
In \citep{chang2006case}, a Genetic Algorithm (GA) is developed for solving $1|r_j|w_jC_j$, and the authors test it on instances with up to 60 jobs.
To measure solution quality, they compute an ``average percent error,” defined as the difference between a given solution and the best solution obtained by a set of  Genetic algorithms. In most instances, after thousands of generations, the algorithm is able to find a solution whose average percent error is below 1\%. However, this metric can only reflect the relative performance between the algorithms considered, and cannot reflect how close the obtained solution is to the true optimum.

For problem $1|r_j|\sum w_jT_j$, only quite small instances have been solved in the literature. \citet{akturk2001new} develops a dominance rule for $1|r_j|\sum w_jT_j$ and integrates it with a total of 11 heuristics to solve instances with 50, 100, and 150 jobs, achieving improvements after applying the rule. However, the optimality gaps of the obtained solutions are not reported. \citet{chou2005heuristic} proposes a heuristic scheduling algorithm, solving instances with relatively small sizes (no more than 20 jobs). Additionally, \citet{jouglet2008dominance} develops dominance-based heuristics for both $1||\sum w_jC_j$ and $1|r_j|\sum w_jC_j$, but only solves small instances with no more than 100 jobs.

In \citep{janiak2009scheduling}, it is shown that $1||w_jC_j^{a_j}$ is NP-hard. The study evaluates five heuristic algorithms on instances ranging from 10 to 500 jobs. For instances with up to 25 jobs, optimal solutions are attainable using a branch-and-bound approach, and the heuristics can find solutions within 1\% of the optimal in just 20–30 milliseconds. However, for larger instances with 50 to 500 jobs, solution qualities are assessed relative to the best-known solutions, and the optimality gaps remain undetermined.

\section{A Unified Machine Scheduling Neural Network}
\label{sec:UMSNN}

In this section, we propose a unified deep neural network (DNN) approach to solve single-machine min-sum scheduling problems. We novelly define the inputs to the DNN using the raw data including job processing times, release dates, and the starting costs in the time-indexed formulation. The neural network then extracts features from the inputs to predict solutions. For ease of presentation, we name this neural network a {\it Unified Machine Scheduling Neural Network} (UMSNN).
The inputs and outputs are defined in subsection \ref{subsec:input_features}. The architecture of UMSNN is presented in subsection \ref{subsec:transformer_encoderdecoder}. Finally, two heuristics for deriving feasible solutions using the outputs of UMSNN are presented in subsection \ref{subsec:list heuristic}.

\subsection{Inputs and outputs}\label{subsec:input_features}

The inputs of the UMSNN should fully describe a given instance. For an instance with a total of $n$ jobs, the inputs are defined as a sequence of $n$ vectors, denoted by $\mathcal{I} = \{I_1, I_2, \dots, I_n\}$, with the $i$-th vector corresponding to job $i$. As presented in Section~\ref{Sec:Introduction}, the time-indexed formulation \eqref{time-index-formulation} is a unified formulation, which incorporates the objective function  into the  starting costs of the jobs. Therefore, we include the starting costs of the jobs in $\mathcal{I}$ to enable our approach to work for all the problems with a min-sum objective function. We note that for a job $j$ with a nonzero release time, its starting cost at every time point $t$ with $t<r_j$, which is calculated as $c_{jt}=z_j(t+p_j)$, is also included in $\mathcal{I}$. Moreover, to learn to generate feasible solutions, the job processing times and the release times, which play critical roles in the constraints of the time-indexed formulation, are also  included in $\mathcal{I}$. 
However, based on our computational experiment, directly utilizing $\{c_{jt}\}_{j\in \mathcal{N},t\in \mathcal{T}}$, $\{p_j\}_{j\in \mathcal{N}}$, and $\{r_j\}_{j\in \mathcal{N}}$ as the inputs does not lead to good results. In the following, an encoding approach is developed to represent these inputs properly. 

For different objective functions, the values of the starting costs may have different orders of magnitude. For example, starting costs in an instance of $1||\sum_{j\in \mathcal{N}} w_j C_j^{a_j}$ with $a_j\ge 1$ may be  significantly larger than those in an instance of $1||\sum_{j\in \mathcal{N}} w_j T_j$. To learn to solve instances with different objective functions by using a single model, we normalize the starting costs $c_{jt}$ as $\Tilde{c}_{jt}=(c_{jt}-\underline{c})/(\Bar{c}-\underline{c})$, where $\underline{c}$ and $\Bar{c}$ are the minimal and maximal value among $\{c_{jt}\}_{j\in \mathcal{N},t\in \mathcal{T}}$, respectively. The normalized starting costs $\{\Tilde{c}_{jt}\}_{j\in \mathcal{N},t\in \mathcal{T}}$ are all within the interval $[0, 1]$. For ease of representation, let $\tilde{c}_j = [\tilde{c}_{j0},\tilde{c}_{j1},..., \tilde{c}_{j,T-1}] \in \mathbb{R}^{1 \times T}$.

The job processing time $p_j$ is encoded as $\tilde{p}_j$, where
\begin{equation}\label{eq:coding_pros}
    \Tilde{p}_j=[\underbrace{1, 1, ..., 1}_{p_j}, 0, ..., 0].
\end{equation}
The length of $\Tilde{p}_j$ should be at least equal to the maximum possible processing time. 
The motivation for encoding the job processing times in this manner, rather than using widely used  one-hot encoding, is to implicitly capture the similarities between jobs. Specifically, if the processing times of two jobs are close, then the two jobs are similar in terms of their processing times. The above encoding scheme (\ref{eq:coding_pros}) reflects this very well because with this encoding scheme, for any two similar jobs $j_1$ and $j_2$, their coded processing times $\Tilde{p}_{j_1}$ and $\Tilde{p}_{j_2}$ would have a lot of overlap.
Moreover, as will be presented in Section~\ref{subsec:special_instances}, the instances for supervised training are specially designed, where job processing times share a common divisor. Therefore, during training, the neural network may only encounter jobs with certain processing times. Through the similarities between the jobs, as represented in the encoded job processing times, the model is expected to be able to predict solutions for  generic instances based on the knowledge learned from special instances.

Similarly, for $j\in\mathcal{N}$, the release date $r_j$ is encoded as 
$\tilde{r}_j = [\tilde{r}_{j0},\tilde{r}_{j1},..., \tilde{r}_{j,T-1}]$, where 
\begin{equation*}
\begin{aligned}
    \Tilde{r}_{jt}=\begin{cases} 
\ 1, & \text{if } t \in \{0, 1, ..., r_j-1\}, \\
\ 0, & \text{otherwise}.
\end{cases}
\end{aligned}
\end{equation*}
For job $j$, $\tilde{r}_j$   has the same length as the planning horizon, and the values associated with the time points before the release date are set to one to indicate the infeasible starting time. For an instance without release dates, all the elements in $\{\Tilde{r}_{jt}\}_{j\in \mathcal{N},t\in \mathcal{T}}$ are set as zero.

The encoded starting costs, processing times, and release dates, are used as the inputs, i.e., $\mathcal{I}=\{I_1, I_2, ..., I_n\}$ with $I_j=(\Tilde{c}_{j},\Tilde{p}_j,\Tilde{r}_j), \forall j \in \mathcal{N}$. To make the resulting neural network work for all instances of all the problems considered, the dimension of each input element needs to be fixed at a large enough number a priori. Thus, we fix the length of the planning horizon $T$ and the maximum possible processing time of the jobs at large enough values beforehand.
This ensures that our UMSNN, which is described  in subsection \ref{subsec:transformer_encoderdecoder}, can solve any instances with any number of jobs and any distribution of job processing times, as long as the length of the planning horizon of the instance and the maximum job processing time in the instance do not exceed the pre-specified values, respectively. 

In the time-indexed formulation \eqref{time-index-formulation}, there are $T$ decision variables $\{x_{jt}\}_{t\in \mathcal{T}}$ for each job $j$. If in our UMSNN, an output of length $T$ is used for each job, there would be too many parameters to be learned such that it would be too time consuming to train the model. To reduce the number of learnable parameters needed, we propose to use time windows (a group of consecutive time points) instead of individual time points as the basic time units associated with the output. Specifically, we divide the planning horizon $\mathcal{T}$ into $\gamma$ time windows  $\mathcal{K}=\{\{0,..., \eta-1\}, \{\eta,..., 2\eta-1\}, ..., \{(\gamma-1)\eta,..., T-1\}\}$, where $\gamma=\lceil \frac{T}{\eta} \rceil$. The output for each job $j$ is a $\gamma$ dimensional vector $o_j = [o_{j0}, o_{j1},..., o_{j,\gamma-1}]\in \mathbb{R}^{1 \times \gamma}$, representing an assignment of the job to the $\gamma$  time windows,  where $o_{jt}$ is the probability that job $j$ is assigned to start at a time point in  the $t$th time window. 

\subsection{Architecture of UMSNN}
\label{subsec:transformer_encoderdecoder}

\begin{figure}[!t]
\begin{center}
        \includegraphics[width=1\linewidth]{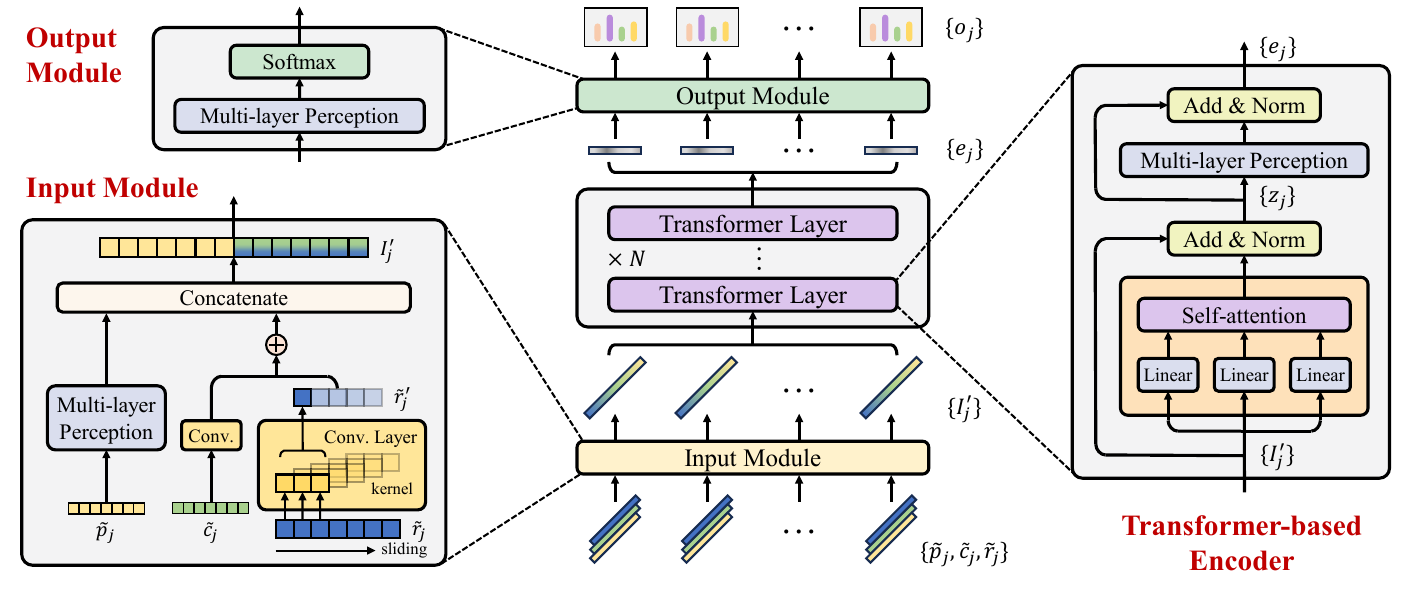}
    \end{center}
    \caption{Architecture of the Unified Machine Scheduling Neural Network}
    \label{Fig:transformer_encoder_decoder}
\end{figure}

In this subsection, we present the architecture of our UMSNN. As presented in the previous subsection, the inputs of the UMSNN are defined based on the time-indexed formulation to make our approach work for various objective functions. However, the inputs defined in this way are raw and unstructured, with a large dimensionality. This means that there can be a vast number of hidden features within the inputs. There  can also be a lot of noise and irrelevant information in the inputs. 
Consequently, it could be difficult for the ML model to generalize, leading to overfitting. To overcome the difficulty, our idea is to leverage the capability of deep neural networks (DNNs) to automatically extract features from the inputs. Specifically,  we use a DNN, which has multiple layers of neurons.  The features within the inputs are progressively extracted, with the low-level features extracted by the earlier layers and complex, high-level features captured by deeper layers. The architecture of the UMSNN is depicted in Figure \ref{Fig:transformer_encoder_decoder}, which has three parts, an input module, a transformer-based encoder, and an output module. 

\textbf{Input module:} In the inputs, $\Tilde{c}_{j}$ and $\tilde{r}_j$ have the same length as the planning horizon (i.e.,  $T$), which can be very large. An input module is developed to reduce the dimensions of $\{I_j\}_{j\in\mathcal{N}}$, while extracting the critical features for predicting solutions.

The input module is a job-wise operation, i.e., it separately processes $I_j$ for each $j$. 
We use multiple convolutional neural network (CNN) layers to extract features from the starting costs $\{\Tilde{c}_{jt}\}_{t\in\mathcal{T}}$. The motivation behind is to exploit the fact that the starting costs among different jobs and instances share similar features, especially for the same objective function. 
To illustrate, take  $1||\sum w_{j}T_{j}$ as an example. In any instance of this problem, the starting costs for jobs  follow a piecewise linear pattern: they start at zero, remain zero up to a certain point, and then increase linearly beyond that. The location of this turning point, which varies among different jobs and instances, is determined by the job processing time and due date, and is critical for predicting solutions. 
A CNN layer utilizes kernels (also known as filters) with learnable parameters that slide over a given input feature vector and perform a multiplying operation at each position. This mechanism enables the layer to have translation-invariant characteristics \citep{lenet1998}. Patterns and features from the input can be identified and extracted irrespective of variations in position. Moreover, since the elements at different positions are processed by the same kernels, the layer has a reasonable number of learnable parameters, even though the input is of high-dimensional. We stack multiple CNN layers and activation functions to extract the features from the starting costs $\{\Tilde{c}_{jt}\}_{t\in\mathcal{T}}$. 

Denote the output vector as $\Tilde{c}_j^{\prime} \in \mathbb{R}^{1\times d^{c}}$ with a dimension of $d^{c}$. The value of $d^{c}$ should be  much smaller than $T$, e.g., in our computational tests (see Section \ref{Sec:Numerical_Results}), $T$ is either 53000 or 68000 whereas $d^{c}$ is only 1024 in both examples. Similarly, for each job $j$, the encoded release date vector $\Tilde{r}_{j}$ is processed by the CNN layers. Denote the output vector as $\Tilde{r}_j^{\prime} \in \mathbb{R}^{1\times d^{r}}$.  By setting the dimensions and stride of the convolution kernel \citep{CNNserveyZeWenLi}, we let the dimension $d^{r}=d^{c}$, which is thus much smaller than $T$. 
The job processing times $\Tilde{p}_j$ are processed by a multilayer perceptron, and the output vector is denoted as $\Tilde{p}_j^{\prime} \in \mathbb{R}^{1\times d^{p}}$, where $d^{p}$ is the dimension. Since the amount of information contained in the processing times $\Tilde{p}_j$ is less than that in the starting costs, the value of $d^{p}$  is set to a value smaller than $d^{c}$. For example, it is set as 256 in Section~\ref{Sec:Numerical_Results}.

The output of the input module is generated by adding $\Tilde{c}_j^{\prime}$ and $\Tilde{r}_j^{\prime}$ and concatenating with $\Tilde{p}_j^{\prime}$, i.e., 
$I^{\prime}_j=[\Tilde{p}_j^{\prime}, \Tilde{c}_j^{\prime}+\Tilde{r}_j^{\prime}] \in \mathbb{R}^{1 \times d^{e}}$, where $d^{e}=d^{p}+d^{c}$.
Again, the input module operates on a job-wise basis, meaning that the input for each job is independently processed by a single input module. Therefore, the module is flexible and works for instances containing a varying number of jobs.

\textbf{Transformer-based Encoder:}
A transformer-based encoder is developed to capture the relations among different jobs. The key component is a self-attention mechanism, which first generates three vectors, called query, key, and value, for each job $j$, as follows:
\begin{equation}
    q_j=I^{\prime}_jW^{Q}, k_j=I^{\prime}_jW^{K}, v_j=I^{\prime}_jW^{V},
\end{equation}
where $W^{Q}, W^{K}, W^{V} \in \mathbb{R}^{d^{e} \times d^{e}}$ are learnable parameters. Attention scores between jobs are then computed. Taking the first job as an example, the attention scores quantifying the relationship between it and all the other jobs are calculated as: 
\begin{equation}
    a_{1j}=\frac{q_1k_j^{T}}{\sqrt{d^{e}}}, \forall j\in \mathcal{N}.
\end{equation}
After normalizing the sum of scores to 1 through a softmax function, the feature vector of the first job is recalculated as $z_1$ by aggregating the features from all the other jobs, as follows:
\begin{equation}
    z_{1}=\sum_{j \in \mathcal{N}}{\frac{\exp(a_{1j})}{\sum_{j' \in \mathcal{N}}\exp(a_{1j'})} \cdot v_j}.
\end{equation}
Vector $z_{1}$ incorporates critical information for predicting a solution for the first job. The above attention operation is applied to every job, which then generates $\{z_j\}_{j\in \mathcal{N}}$.

Two commonly used techniques known as the residual connection and the layer normalization techniques \citep{ba2016layer,vaswani2017attention} are used to prevent gradient vanishing or gradient exploding in the transformer architecture.
The two techniques are depicted as ``Add \& Norm" in Figure \ref{Fig:transformer_encoder_decoder}. After that, a feed forward layer is used to further refine and extract features. Then, the residual connection and the layer normalization are used again. Denote the outputs of the encoder as $\{e_j\}_{j\in\mathcal{N}}$. 

\textbf{Output module:}
Based on $\{e_j\}_{j\in\mathcal{N}}$, an output module generates predictions. For job $j$, $e_j$ is first processed through a multilayer perceptron, with the resulting output denoted as $e_j^{\prime} \in \mathbb{R}^{\gamma \times 1}$. Subsequently, an output vector $o_j \in \mathbb{R}^{\gamma\times 1}$ is produced using a softmax function, where the $k$th element is calculated as: 
\begin{equation}\label{UMSNN_output}
    o_{jk}=\frac{\exp(e^{\prime}_{jk})}{\sum_{\kappa=1}^{\gamma}\exp(e^{\prime}_{j\kappa})}, \forall ~k  = 1,2, ..., \gamma.
\end{equation}
The elements of vector $o_j$ lie between zero and one, summing up to one. Therefore, $o_j$ is treated as a probability vector, with its $k$th element, $o_{jk}$, interpreted as the probability that the optimal starting time of job $j$ falls within the $k$th time window (i.e., $\{k\eta, ..., (k+1)\eta\}$).

For ease of the notation, we denote the learnable parameters in the UMSNN collectively as $\theta\in \mathbb{R}^{\zeta}$, where $\zeta$ is the corresponding dimension, and we denote the UMSNN itself as $f^{U}_{\theta}(\cdot)$, which maps a given input $I$ to an output $o$.

\subsection{Heuristics for feasibility}
\label{subsec:list heuristic}
The outputs of UMSNN are probability vectors $\{o_j\}_{j\in\mathcal{N}}$, with $o_j$ associated with job $j$ and has a length of $\gamma$. In this subsection, two simple heuristics are developed to generate a feasible starting time for each job based on the probability vectors. Both approaches first sequence the jobs in a certain order, and then calculate the starting times of the jobs accordingly.

The first approach finds the time window with the largest probability for each job $j\in \mathcal{N}$, i.e.,
\begin{equation}
    \omega_j=\argmax_{k\in\mathcal{K}} \{o_{jk}\}, 
\end{equation}
and sequences the jobs  by the ascending order of $\{\omega_j\}_{j\in\mathcal{N}}$. If multiple jobs fall within the same time window, they are sequenced by the descending order of their weights. Since only the time window with the largest probability is considered, we call this approach a ``greedy" approach. 

It is possible that for a job, a time window that does not have the largest probability has a high quality. Therefore, our second approach may assign any time window $k$ with a positive $o_{jk}$ value to  a given job $j$. More specifically, it is a sampling-based approach that chooses  the $k$th time window for job $j$ with probability $o_{jk}$.
Through sampling a time window for each job, $\{\omega_j\}_{j\in\mathcal{N}}$ is obtained, and the jobs are sequenced following the ascending order of the time windows as in the first approach. This approach is called a ``sampling" approach.

In both approaches, given the sequence of the jobs, feasible starting times for the jobs can be easily calculated. For the first job, the starting time is set as its release date (which is zero if the problem does not consider release times). For the successive jobs, the starting time equals the greater value between its release time and the completion time of its predecessor.

\section{Supervised Training Using Special Instances}
\label{sec:generation_of_trainingset}

To generate high-quality solutions for online instances, the learnable parameters of the UMSNN, as defined in the previous section, need to be trained using a large variety of instances  generated from the same distribution as the online instances to be solved. For a given instance, the objective value of a solution generated by the UMSNN, i.e., the total cost of the jobs,  can measure the quality of the solution. However, numerical tests suggest that, if the objective values are  used in the performance measure in offline training, the performance of the UMSNN for online instances would not be satisfactory. This is because the objective value of a solution is only an aggregated piece of information based on the outputs $\{o_{jk}\}_{j\in \mathcal{N}, k\in \mathcal{K}}$ from the UMSNN corresponding to this solution, and hence contains far less information than the entire set of the outputs.

In this paper, we use the entire set of outputs of the UMSNN, $\{o_{jk}\}_{j\in \mathcal{N}, k\in \mathcal{K}}$,  to define a performance measure to be optimized in the offline training stage. Specifically, for a given training instance (or input vector) $\mathcal{I}$, suppose that the outputs of the UMSNN are $\{o_{jk}(\mathcal{I})\}_{j\in\mathcal{N}, k\in\mathcal{K}}$, where   $o_{jk}(\mathcal{I})$ represents the probability that the starting time of job $j$ falls within the $k$th time window in $\mathcal{K}$, and suppose that $\{o^{\star}_{jk}(\mathcal{I})\}_{j\in\mathcal{N},k\in\mathcal{K}}$ are the labels denoting the optimal solution such that  they follow a one-hot format, i.e., $o^{\star}_{jk}(\mathcal{I})=1$ if job $j$ starts in the $k$th time window in the optimal solution, and 0 otherwise. We use the cross-entropy loss \citep{goodfellow2016deep} to measure the error of the outputs from the UMSNN. 
The cross-entropy loss associated with job $j$ is
$
    \sum_{k\in\mathcal{K}}o_{jk}^{\star}(\mathcal{I})\ln {o_{jk}(\mathcal{I})}.
$
Suppose that the optimal starting time of job $j$ falls within the $k^{\star}$th  time window. Then, the loss becomes $ -\ln {o_{jk}(\mathcal{I})}$. 
Minimizing it is equivalent to maximizing the value of $o_{jk}^{\star}(\mathcal{I})$, which is the probability that the optimal time window is correctly predicted. The cross-entropy loss of instance $\mathcal{I}$ is the average loss over all the jobs, i.e.,
\begin{equation}
\label{eq:cross_entropy_loss}
    \mathcal{L}(o^{\star}(\mathcal{I}),o(\mathcal{I}))=-\frac{1}{n}\sum_{j\in\mathcal{N}}\sum_{k\in\mathcal{K}}o_{jk}^{\star}(\mathcal{I})\ln {o_{jk}(\mathcal{I})}.
\end{equation}
With the loss function defined above, the learnable parameters $\theta$ are optimized by solving the following optimization problem:
\begin{equation}
\label{goal_of_offline}
  \min_{\theta} \sum_{\mathcal{I}\in \mathcal{D}} \mathcal{L}(o^{\star}(\mathcal{I}),o(\mathcal{I})),
\end{equation}
where $\mathcal{D}$ is the training instance set. This problem can be solved using the stochastic gradient descent algorithm \citep{amari1993backpropagation} where the gradients of the loss functions with respect to $\theta$ can be calculated by the backpropagation algorithm discussed in Appendix~\ref{sec:back_propagation}. 

The biggest challenge in supervised training of our UMSNN is generating labels (i.e., optimal time windows to start the jobs) for training instances, especially when the training instances have large sizes. 
This is because the scheduling problems considered in this paper are strongly NP-hard. To overcome this difficulty, we develop an approach to generate a set of special instances for which optimal solutions can be found quickly. This is presented in subsection \ref{subsec:special_instances}. In subsection \ref{subsec:data_augment}, a data augmentation approach is developed to further enrich the training set without extra computational effort.

\subsection{Generation of special instances} 
\label{subsec:special_instances}

As discussed in Section~\ref{Sec:Introduction}, the time-indexed formulation \eqref{time-index-formulation} often generates a tight LP relaxation bound. Thus, for instances with a relatively small size in terms of the number of binary variables $n\cdot T$, the time-indexed formulation can be solved to optimality quickly. However, it is not practical to use this formulation to solve large instances directly.
This motivates us to develop an approach to generate  large instances based on much smaller instances such that (i) the time-indexed formulations of the smaller instances can be solved to optimality quickly, and (ii) the optimal solutions of the smaller instances can be easily converted to optimal solutions of the  larger instances. To illustrate the idea of our approach, consider a small instance of $1||\sum T_j$ with the values of the job processing times and due dates being small, e.g., they are all in the set $\{1, 2, ..., 10\}$. 
 We convert this instance to a much larger instance which contains the same number of jobs but with the processing time $p_j^{\prime}$ and due time $d_j^{\prime}$ of each job $j$ being $\alpha$ times of the original values $p_j$ and $d_j$, respectively, where $\alpha$ is a positive integer. Thus, in the larger instance, job processing times and due dates all belong to the set $\{\alpha, 2\alpha, ..., 10\alpha\}$.  Clearly, the size of the time-indexed formulation for the smaller instance is only $1/\alpha$ of that for the larger instance. As proved later, given the optimal starting times of jobs of the smaller instance, denoted as $\{S_1, ..., S_n\}$, an optimal solution to the larger instance can be easily constructed as $\{\alpha  S_1, ..., \alpha  S_n\}$. Since the parameters of the large instances have some special characteristics (e.g., have a common divisor $\alpha$), we refer to these instances as special instances.

In the following, we summarize our idea as a theorem, based on which special larger instances can be constructed from randomly generated smaller instances. To make our idea as general as possible, the following theorem is based on the time-indexed formulation \eqref{time-index-formulation}, where different objective functions are unified by the starting costs of the jobs.

\begin{theorem} 
\label{thm1}
For any given single-machine scheduling problem with a min-sum non-decreasing objective function, suppose that we are given an instance consisting of the following parameters, among others: set of jobs $\mathcal{N}=\{1, ..., n\}$,   planning horizon  $\mathcal{T}=\{0, 1, \ldots, T-1\}$,  job processing times  $\{p_j\}_{j\in\mathcal{N}}$, job release times  $\{r_j\}_{j\in\mathcal{N}}$, and  starting costs of the jobs  $\{c_{jt}\}_{j\in\mathcal{N}, t\in\mathcal{T}}$. Suppose that this instance is solved optimally with the optimal starting times of jobs as $\{S_1, \ldots, S_n\}$.
Given any positive integer $\alpha$, construct a larger instance, where there is the same number of jobs as in the smaller instance, but the  planning horizon $\mathcal{T}^{\prime}$, and the processing time $p_j^{\prime}$ and release time $r_j^{\prime}$ of each job $j\in \mathcal{N}$ are all enlarged to $\alpha$ times the corresponding values in the smaller instance, i.e., $\mathcal{T}^{\prime} = \{0, 1, \ldots, \alpha (T - 1)\}$, $p_j^{\prime}=\alpha p_j$ and $r_j^{\prime}=\alpha r_j$ for $j\in \mathcal{N}$. In addition, the values of some other parameters (e.g., due dates) could also be adjusted for the larger instance.
If there exists a scalar $\beta$ such that the resulting starting costs in the larger instance $\{c_{jt}^{\prime}\}_{j\in\mathcal{N}, t\in\mathcal{T}^{\prime}}$ satisfy:
\begin{equation}\label{eq_for_specialinstance}
    c^{\prime}_{j,\alpha t}= \beta \cdot c_{jt}, ~~~~~ \mbox{for each $j\in\mathcal{N},t\in \mathcal{T}$},
\end{equation}
then using $\{\alpha S_1, \ldots, \alpha S_n\}$ as starting times of the jobs gives an optimal solution to the larger instance.
\end{theorem}
\proof{Proof}
We prove the theorem by contradiction. Since the objective functions are non-decreasing in the starting times of the jobs, and the processing times and release times of the jobs in the larger instance share a common divisor of $\alpha$, the larger instance has an optimal solution where the starting times of the jobs have a common divisor $\alpha$.
Suppose that 
$\{\alpha S_j\}_{j \in \mathcal{N}}$ is not optimal to the larger instance, and its true optimal solution is $\{\alpha \tilde{S}_j\}_{j \in \mathcal{N}}$. This implies the following inequality
\begin{equation}\label{proof:assumption}
    \sum_{j\in \mathcal{N}}c^{\prime}_{j,\alpha  \tilde{S}_j}<\sum_{j\in \mathcal{N}}c^{\prime}_{j,\alpha  S_j}.
\end{equation}
This, along with \eqref{eq_for_specialinstance}, implies that
\begin{equation}
    \sum_{j\in \mathcal{N}}c_{j,\tilde{S}_j}<\sum_{j\in \mathcal{N}}c_{j,S_j},
\end{equation}
which indicates that  $\{\tilde{S}_j\}_{j \in \mathcal{N}}$ is a better solution to the given instance than  $\{S_j\}_{j \in \mathcal{N}}$,  leading to a contradiction. 
\endproof

We note that when generating a larger instance from a given smaller instance by applying Theorem~\ref{thm1}, in order to satisfy \eqref{eq_for_specialinstance}, other job parameters, besides  job processing times and release times,  may also need to be enlarged accordingly.  For example, for $1||\sum w_jT_j$, we also need to enlarge job due dates to $\alpha$ times the corresponding values in the smaller instance.  For some problems, job weights may also need to be adjusted to satisfy \eqref{eq_for_specialinstance}. For example, for
 $1||w_jC_j^{a_j}$, to satisfy \eqref{eq_for_specialinstance}, the weight $w_j$ needs to be reduced to $\alpha^{-a_j}$ times the original values. 
Table \ref{tab:Approach to generating special instances} shows how larger instances can be generated from given smaller instances for the nine  problems to be used in our computational experiment. 
\begin{table}[ht]
    \centering 
    \caption{Approach to generating larger instances} 
    \label{tab:Approach to generating special instances}
    \begin{adjustbox}{max width=\linewidth}
    \begin{tabular}{ccccc}  
        \toprule 
    \multirow{2}*{Problem}  & \multicolumn{2}{c}{Given Small Instances}& \multicolumn{2}{c}{Generated Larger Instances}  \\
    \cmidrule(lr){2-3} \cmidrule(lr){4-5}
    & \makecell{Parameters} &\makecell{Solution} & \makecell{Parameters} &\makecell{Solution}\\
     \midrule 
        $1|r_j|\sum w_j C_j$ &$\{p_j,r_j, w_j\}_{j\in\mathcal{N}}$& $\{S_j\}_{j\in\mathcal{N}}$ &$\{\alpha p_j, \alpha  r_j, w_j\}_{j\in\mathcal{N}}$& $\{\alpha  S_j\}_{j\in\mathcal{N}}$\\
        \midrule 
        $1||\sum w_j C_j^{a_j}$ &$\{p_j,a_j,w_j\}_{j\in\mathcal{N}}$& $\{S_j\}_{j\in\mathcal{N}}$ &$\{\alpha p_j, a_j,  w_j/\alpha^{a_j}\}_{j\in\mathcal{N}}$& $\{\alpha  S_j\}_{j\in\mathcal{N}}$\\
        $1|r_j|\sum w_j C_j^{a_j}$ &$\{p_j,r_j, a_j, w_j\}_{j\in\mathcal{N}}$& $\{S_j\}_{j\in\mathcal{N}}$ &$\{\alpha p_j, \alpha  r_j, a_j, w_j/\alpha^{a_j}\}_{j\in\mathcal{N}}$& $\{\alpha  S_j\}_{j\in\mathcal{N}}$\\
        \midrule
        $1||\sum w_j T_j$ &$\{p_j,d_j,w_j\}_{j\in\mathcal{N}}$& $\{S_j\}_{j\in\mathcal{N}}$ &$\{\alpha p_j,  \alpha  d_j, w_j\}_{j\in\mathcal{N}}$& $\{\alpha  S_j\}_{j\in\mathcal{N}}$\\
        $1|r_j|\sum w_j T_j$ &$\{p_j,r_j,d_j,w_j\}_{j\in\mathcal{N}}$& $\{S_j\}_{j\in\mathcal{N}}$ &$\{\alpha p_j,\alpha  r_j,  \alpha  d_j, w_j\}_{j\in\mathcal{N}}$& $\{\alpha  S_j\}_{j\in\mathcal{N}}$\\
        \midrule
        \makecell{$1||\rho\sum  w_{1j}T_j+(1-\rho) \sum w_{2j}C_j$} &$\{p_j,d_j,w_{1j},w_{2j}\}_{j\in\mathcal{N}}$& $\{S_j\}_{j\in\mathcal{N}}$ &$\{\alpha p_j, \alpha  d_j, w_{1j},w_{2j}\}_{j\in\mathcal{N}}$& $\{\alpha  S_j\}_{j\in\mathcal{N}}$\\
        \makecell{$1|r_j| \rho \sum w_{1j}T_j+(1-\rho)\sum w_{2j}C_j$} &$\{p_j,r_j, d_j, w_{1j},w_{2j}\}_{j\in\mathcal{N}}$& $\{S_j\}_{j\in\mathcal{N}}$ &$\{\alpha p_j, \alpha  r_j, \alpha  d_j, w_{1j},w_{2j}\}_{j\in\mathcal{N}}$& $\{\alpha  S_j\}_{j\in\mathcal{N}}$\\
        \midrule
        \makecell{$1||\rho\sum w^{A}_jT^{A}_j$$+(1-\rho)\sum w^{B}_jC^{B}_j$} &\makecell{$\{p_j,w_j\}_{j\in\mathcal{N}_1}$\\$\{p_j,d_j, w_j\}_{j\in\mathcal{N}_2}$}&$\{S_j\}_{j\in\mathcal{N}}$& \makecell{$\{\alpha p_j,w_j\}_{j\in\mathcal{N}_1}$\\$\{\alpha p_j,\alpha  d_j, w_j\}_{j\in\mathcal{N}_2}$}&$\{\alpha  S_j\}_{j\in\mathcal{N}}$\\
        \makecell{$1|r_j|\rho\sum w^{A}_jT^{A}_j$$+(1-\rho)\sum w^{B}_jC^{B}_j$}& \makecell{$\{p_j,r_j,w_j\}_{j\in\mathcal{N}_1}$\\$\{p_j,r_j,d_j,w_j\}_{j\in\mathcal{N}_2}$}&$\{S_j\}_{j\in\mathcal{N}}$& \makecell{$\{\alpha p_j, \alpha  r_j,w_j\}_{j\in\mathcal{N}_1}$\\$\{\alpha p_j,\alpha  r_j,\alpha d_j,w_j\}_{j\in\mathcal{N}_2}$}&$\{\alpha  S_j\}_{j\in\mathcal{N}}$\\
        \bottomrule 
    \end{tabular}
    \end{adjustbox}
\end{table}
Given any problem included in 
Table \ref{tab:Approach to generating special instances}, we can generate a set of special instances and their optimal solutions by (i) first randomly generating a set of small instances with short job processing times (and short job release dates if the problem involves release dates), (ii) then solve each of them using the time-indexed formulation, and (iii) finally, constructing the corresponding larger instances and optimal solutions following Table \ref{tab:Approach to generating special instances}.

\subsection{Data augmentation} 
\label{subsec:data_augment}
To enrich the instances for training, we develop a data augmentation approach  by exploiting the fact that the optimal solution of an instance may remain optimal when some parameters in the instance (e.g., due dates or release dates) are changed slightly.
Our idea is to try to utilize instances that are already generated and solved with optimal solutions to generate new instances that have the same optimal solutions with little extra computational effort by making slight changes to  some parameters in the instances. In the following, we first present a theorem as the basis for our approach. The implementation of this approach is described in Section ~\ref{Sec:Numerical_Results} where we report computational results. 

\begin{theorem}
    For any given single-machine scheduling problem with a min-sum non-decreasing objective function, suppose that we are given an instance consisting of the following parameters, among others: set of jobs $\mathcal{N}=\{1, ..., n\}$,   planning horizon  $\mathcal{T}=\{0, 1, \ldots, T-1\}$,  and job processing times  $\{p_j\}_{j\in\mathcal{N}}$,   job release times  $\{r_j\}_{j\in\mathcal{N}}$ (which could all be 0),  and  starting costs of the jobs  $\{c_{jt}\}_{j\in\mathcal{N}, t\in\mathcal{T}}$. Suppose that this instance is solved optimally with the optimal starting times of jobs as $\{S_1, \ldots, S_n\}$.
    Construct a new instance of the same problem with the same jobs, same planning horizon, and same jobs processing times, but with  modified release times $\{r^{\prime}_{j}\}_{j\in\mathcal{N}}$ such that: 
    \begin{equation}
    r_{j} \leq r^{\prime}_{j}\leq S_j, ~~\forall j \in\mathcal{N},  \label{eq_rj}
    \end{equation}
    and possibly some other parameters (e.g., due dates $d_j$ if they are relevant) also modified.  If the resulting starting costs $\{c^{\prime}_{j,t}\}_{j\in\mathcal{N},t\in\mathcal{T}}$ in the new instance satisfy the following: 
    \begin{equation}  \label{eq_cj}
       c^{\prime}_{j,S_j}= c_{j,S_j}, ~~  c^{\prime}_{j,t}\geq c_{j,t} ~~\forall t\in\mathcal{T}\setminus\{S_j\}, ~\forall j\in\mathcal{N}
    \end{equation}
    then $\{S_j\}_{j\in\mathcal{N}}$ is also an optimal solution to the new instance.    
\end{theorem}

\proof{Proof}
We prove the theorem by contradiction. 
Suppose that for the new instance, 
$\{S_j\}_{j \in \mathcal{N}}$ is not optimal, and its true optimal solution is $\{\tilde{S}_j\}_{j \in \mathcal{N}}$. Thus, $\tilde{S}_j\ge r_j^{\prime}$, and hence by \eqref{eq_rj}, we have $\tilde{S}_j\ge r_j$, for $j \in\mathcal{N}$. This implies that  $\{\tilde{S}_j\}_{j \in \mathcal{N}}$ is a feasible solution to the original instance. Since $\{\tilde{S}_j\}_{j \in \mathcal{N}}$ is optimal for the new instance, we have
\begin{equation}\label{proof:thm2-assumption}
    \sum_{j\in \mathcal{N}}c^{\prime}_{j,\tilde{S}_j}<\sum_{j\in \mathcal{N}}c^{\prime}_{j,S_j}.
\end{equation}
This, along with \eqref{eq_cj}, implies that
\begin{equation}
    \sum_{j\in \mathcal{N}}c_{j,\tilde{S}_j}\leq \sum_{j\in \mathcal{N}}c^{\prime}_{j,\tilde{S}_j}<\sum_{j\in \mathcal{N}}c^{\prime}_{j,S_j} =
    \sum_{j\in \mathcal{N}}c_{j,S_j},
\end{equation}
which indicates that  $\{\tilde{S}_j\}_{j \in \mathcal{N}}$ is a better solution to the original instance than  $\{S_j\}_{j \in \mathcal{N}}$,  leading to a contradiction. 
\endproof

Based on the above theorem, given any instance, new instances with the same optimal solution  can be generated by changing the release times or / and some other parameters such that the starting costs in the new instance satisfy \eqref{eq_cj}. For a problem with a due-date related objective function, e.g., $1|r_j|\sum w_jT_j$, we can generate new instances from a given instance by modifying (i) the release times  following the theorem, and (ii) the value of due date $d_j$ of any job $j$ that is completed before its due date in the optimal solution (i.e., $S_j+p_j< d_j$) for the original instance to a new value $d_j^{\prime}$ between $S_j+p_j$ and $d_j$, because such  changes  satisfy \eqref{eq_cj}.

\section{Online Single-Instance Learning} \label{sec:online_learning}

The offline training approach  developed in the previous section trains the UMSNN model for optimized average performance on a large number of offline instances of one or more given single-machine min-sum scheduling problems. Given an online instance to be solved, it is unlikely that the offline trained parameters are optimal  for this instance. Therefore, in this section, we propose an approach to fine-tune the learnable parameters  whenever a given online instance needs to be solved. Since fine-tuning is based on a single instance only, our approach is called ``single-instance" learning. A key element of the approach is a feasibility surrogate that we develop, which connects the output of the UMSNN to a solution of the underlying scheduling problem.
We first describe our overall approach in subsection \ref{subsec:Motivation and the entire framework} and then describe the feasibility surrogate  in subsection \ref{subsec:feasibility surrogate}. 

\subsection{Approach}
\label{subsec:Motivation and the entire framework}

The UMSNN model $f^U_{\theta}(\cdot)$, once trained (i.e., given the values of the learnable parameters $\theta\in \mathbb{R}^{\zeta}$), can map any given instance, represented by the input vector $\mathcal{I}$, to an output vector $o = f^U_{\theta}(\mathcal{I})$, where  row $o_j$ is the probability vector of the starting time of job $j$ assigned to the $\gamma$ time windows, as defined in Section~\ref{sec:UMSNN}. Now, suppose that the input vector $\mathcal{I}$ is given and fixed, but the learnable parameters $\theta$ are variables. We can then view $f^U_{\theta}(\mathcal{I})$ as a function of $\theta$, denoted as $G_{\mathcal{I}}(\theta)$, that maps the learnable parameters $\theta$  to output $f^U_{\theta}(\mathcal{I})$, i.e., $G_{\mathcal{I}}(\theta) = f^U_{\theta}(\mathcal{I})$.  Furthermore, suppose that we have a feasibility layer $F_{\mathcal{V}}(\cdot)$, which can map any output $o$ of the UMSNN to a feasible solution of the job starting times $S=[S_1, S_2, ..., S_n]\in \mathbb{R}^{n}$ based on the release times and the processing times of the given instance $\mathcal{I}$, denoted by $\mathcal{V}=\{p_j, r_j\}_{j\in\mathcal{N}}$.  Such a feasibility layer is developed in subsection \ref{subsec:feasibility surrogate}. Using $F_{\mathcal{V}}(\cdot)$, we can then define the following optimization problem for the given instance $\mathcal{I}$. 
\begin{equation}
\label{alternative_function}
    \min_{\theta\in \mathbb{R}^{\zeta}} J(\theta)=Z(F_{\mathcal{V}}(G_{\mathcal{I}}(\theta))),
\end{equation}
where $Z(S) = \sum_{j \in \mathcal{N}} z_j(S_j + p_j)$ is the objective value. Problem \eqref{alternative_function} is to find optimal $\theta$, denoted as $\theta^\star$, such that for the given instance $\mathcal{I}$, the objective value $J(\theta^\star)$ is minimum. Therefore, problem \eqref{alternative_function} can be viewed as an alternative formulation of the original machine scheduling problem \eqref{time-index-formulation}, but with continuous decision variables $\theta\in \mathbb{R}^{\zeta}$. 

Nevertheless, in view of the  large number of decision variables $\theta$ and the complex relationship between the objective function and the underlying UMSNN model, directly optimizing the alternative problem \eqref{alternative_function} from scratch is impractical. Instead, we use the trained values of $\theta$, denoted as $\theta^{offline}$,  resulted from the offline training of the UMSNN model as the initial solution of $\theta$ to solve the alternative problem \eqref{alternative_function}. Starting with the initial solution $\theta_0=\theta^{offline}$, our gradient descent  algorithm updates the learnable parameters in the $k+1$st iteration, for $k=0, 1, ...,$ as: $\theta_{k+1}=\theta_{k}-s_{k}\cdot g_{k}$. A straightforward choice for $g_{k}$ is the gradient $\nabla_{\theta} J(\theta_{k})$. However, to derive useful  (i.e., nonzero) gradients, there is a challenge we must overcome. Scheduling problems are discrete optimization problems. Thus, for a given instance $\mathcal{I}$, there are finite number of feasible starting times associated with $F_{\mathcal{V}}(G_{\mathcal{I}}(\theta))$, whereas, the variables $\theta$ are continuous real values. As a result, the gradient of the feasibility layer $F_{\mathcal{V}}(\cdot)$ would be zero at most places, as it maps continuous values to discrete ones. To address this challenge, we design a ``feasibility surrogate" $\tilde{F}_{\mathcal{V}}(\cdot)$ to approximate $F_{\mathcal{V}}(\cdot)$. Specifically, $\tilde{F}_{\mathcal{V}}(\cdot)$ approximates feasible job starting times rather than generating them exactly, such that useful gradients can be derived. This is described in subsection~\ref{subsec:feasibility surrogate}. Accordingly, function $J(\theta)$ is then approximated by $\tilde{J}(\theta)=Z(\tilde{F}_{\mathcal{V}}(G_{\mathcal{I}}(\theta)))$. Consequently, in the $k+1$st iteration, for $k=0, 1, ...,$ the learnable parameters  $\theta$ are updated as:
\begin{equation} \label{online_update_iteration}
   \theta_{k+1}=\theta_{k}-s_{k}\cdot \nabla_{\theta} \tilde{J}(\theta_k).
\end{equation}

We now discuss how the gradient $\nabla_{\theta} \tilde{J}(\theta_k)$ is derived. Given the complex structure of the feasibility surrogate $\tilde{F}_{\mathcal{V}}(\cdot)$, it is impractical to explicitly derive the gradient function based on it, i.e., $\nabla_{o} F_{\mathcal{V}}(o)$. For given $o$, we use computational graphs to establish the dependencies among variables and compute derivatives numerically using the chain rule. The gradient calculated is further back-propagated to the UMSNN to calculate the gradient of $\tilde{J}(\theta)$ at point $\theta_k$, i.e., $\nabla_{\theta} \tilde{J}(\theta_k)$.

Next, we discuss how  the stepsize $s_k$ in \eqref{online_update_iteration} is selected. Commonly used stepsizing rules, e.g., the constant stepsize  and the decaying stepsize, are easy to implement, but usually need to be fine-tuned for given problems. For the single-machine scheduling problems we study, different objective functions may have  different magnitudes of values. To have a good performance for different problems without the need to fine tune the stepsize, we use the Polyak stepsize \citep{polyak1969minimization}:
\begin{equation}
\label{online_polyak}
    s_k=\frac{\tilde{J}(\theta_k)-\Bar{J}_{k}}{||\nabla_{\theta} \tilde{J}(\theta_k)||_2^2},
\end{equation}
where $\Bar{J}_{k}$ is a dynamically adjusted target value, and is an estimation of the optimal objective value. Such a stepsize is adaptive  since it is calculated based on the gap between the current objective value $\tilde{J}(\theta_k)$ and the target $\Bar{J}_{k}$. Therefore, the impact of the objective function type on the updating of $\theta$ is reduced. Similar to \citet{goffin1999convergence,nedic2001incremental}, we set the target at iteration $k$ as the best objective value found so far minus a bias $\lambda_k$ as 
\begin{equation}
 \Bar{J}_{k}=\min_{\kappa<k} \tilde{J}(\theta_\kappa) - \lambda_k.
\end{equation}
If the objective value $\tilde{J}(\theta_k)$ does not decrease for successive $B$ iterations, $\lambda_k$ is halved.

Upon the termination of online learning, let $\theta^{\mathrm{online}}$ denote the  solution of $\theta$ found. Feasible job starting times with an improved quality are calculated based on  
$\theta^{\mathrm{online}}$ as
$F_{\mathcal{V}}(G_{\mathcal{I}}(\theta^{\mathrm{online}}))$.

\subsection{Feasibility Surrogate}
\label{subsec:feasibility surrogate}

In subsection \ref{subsec:list heuristic}, we give two heuristic procedures to  generate feasible job starting times $\{S_j\}_{j\in\mathcal{N}}$ based on the output of the UMSNN. These procedures involves multiple non-differentiable operations, such as $\argmax(\cdot)$, the sampling operation, and the job sorting operation. These non-differentiable operations are hard to approximate by differentiable operations. In this subsection, we use a different heuristic to  generate feasible job starting times based on the output of the UMSNN, and use this as the feasibility layer $F_{\mathcal{V}}(\cdot)$ in our single-instance learning approach described in subsection~\ref{subsec:Motivation and the entire framework}. Fewer non-differentiable operations are involved in this heuristic than the two given in subsection \ref{subsec:list heuristic} such that $F_{\mathcal{V}}(\cdot)$ is less difficult to approximate. In the following, we first describe this heuristic, and then briefly explain how the non-differentiable operations involved in the heuristic can be approximated by differentiable operations. 

\textbf{Feasibility layer $F_{\mathcal{V}}(\cdot)$: } The input to $F_{\mathcal{V}}(\cdot)$ is a matrix $o\in\mathbb{R}^{n\times \gamma}$, where $o_{jk}$ is the  probability that the starting time of job $j$ falls within  the $k$th time window. Giving the input $o$, the function $F_{\mathcal{V}}(o)$ consists of two mathematical operations. The first operation calculates the weighted time window value for each $j\in\mathcal{N}$,  as
\begin{equation}
    \tilde{\omega}_j=\sum_{k\in\mathcal{K}}k\cdot o_{jk},
\end{equation}
The second operation  generates feasible job starting times $\{S_j\}_{j\in\mathcal{N}}$, by following the ascending order of the weighted time window values of the jobs $\{\tilde{\omega}_j\}_{j\in\mathcal{N}}$. We assume that the weighted time window values of the jobs are all different. When multiple jobs have the same weighted time window value,   sufficiently small perturbations are added to differentiate them. To express the second operation mathematically, we establish the relationship between $\{\tilde{\omega}_j\}_{j\in\mathcal{N}}$ and the feasible job starting times $\{C_j\}_{j\in\mathcal{N}}$ by using step functions. Specifically, a step function $u(x)$, which is defined as 1 if $x>0$ and 0 otherwise, is used to indicate the relative orders of two jobs. For $j_{1}, j_{2} \in\mathcal{N}$, if $u(\tilde{\omega}_{j_2}-\tilde{\omega}_{j_1})=1$, then $\tilde{\omega}_{j_1}<\tilde{\omega}_{j_2}$ and hence job $j_1$ is ordered before $j_2$, and if $u(\tilde{\omega}_{j_2}-\tilde{\omega}_{j_1})=0$, then $\tilde{\omega}_{j_1}>\tilde{\omega}_{j_2}$ and hence job $j_1$ is ordered after $j_2$. 

To derive job starting times,  we first consider the case where the jobs in the given instance do not have release times. Since in this case, no ideal time should be inserted in an optimal schedule, a job starts immediately after its predecessor is completed. Thus, for each job $j$, its starting time $\tilde{S}_j$ can be calculated by summing up the processing times of all its predecessors, i.e., 
\begin{equation}
\label{eq:feasible_startingtime_without_release}
\tilde{S}_j=\sum_{j^{\prime}\in\mathcal{N}\setminus \{j\}} u(\tilde{\omega}_{j}-\tilde{\omega}_{j^{\prime}})p_{j^{\prime}}, ~~\forall j\in\mathcal{N}.
\end{equation}
A feasible objective value is then calculated as $\sum_{j\in\mathcal{N}}z_j(\tilde{S}_j+p_j)$. 

We now consider the case where jobs have nonzero release times. In this case, a job may not be able to start right after the completion of its predecessor because its release time and the release times of the jobs scheduled before it may push it backward for some time. To derive the starting time of each job in this case, we first derive a formula to calculate the amount of time by which each job needs to be pushed backward in the given job sequence following the ascending order of their weighted time window values $\{\tilde{\omega}_j\}_{j\in\mathcal{N}}$. Without considering job release times, the starting times of the jobs are calculated by \eqref{eq:feasible_startingtime_without_release}. Define $\Delta S_{j}$ to be the amount of time by which  job $j\in \mathcal{N}$ needs to be pushed backward after the job release times are considered. Thus, job $j$'s actual starting time becomes $S_j = \tilde{S}_j + \Delta S_{j}$, for $j\in \mathcal{N}$.

In the following we show how $\Delta S_{j}$ can be calculated. For ease of presentation, we denote the job in the $k$th position of the given sequence as $[k]$, for $k=1, \ldots, n$. For the first job $[1]$, either $\tilde{S}_{[1]}\leq r_{[1]}$ or $\tilde{S}_{[1]}> r_{[1]}$. In the former case, job $[1]$'s starting time should be pushed backward for $r_{[1]}-\tilde{S}_{[1]}$ time slots, and in the latter case, its starting time does not need to be pushed backward. Therefore, $\Delta S_{[1]}$ should be calculated as
\begin{equation}
    \Delta S_{[1]}=\max\{0,r_{[1]}-\tilde{S}_{[1]}\}.
\end{equation}
Now, consider any job $[k]$ in the given sequence, for $k\geq 2$. The fact that job $[k-1]$, which is sequenced immediately before job $[k]$, is pushed backward for  $\Delta S_{[k-1]}$ time units, job $[k]$ must be pushed backward for at least this much time. In the meantime, due to its release time, job $[k]$ needs to be pushed backward for at least $\max\{0, r_{[k]} - \tilde{S}_{[k]}\}$ time units. Thus, $\Delta S_{[k]}=\max\{\Delta S_{[k-1]}, r_{[k]}-\tilde{S}_{[k]}\}$. By recursion, we have
\begin{equation} \label{eq:delta_S[k]}
    \Delta S_{[k]}=\max\{0, (r_{[1]}-\tilde{S}_{[1]}), ..., (r_{[k]}-\tilde{S}_{[k]})\}.
\end{equation}
By using the step function $u(\cdot)$ defined earlier, and using the original job index $j$ instead of their positional index $[j]$, \eqref{eq:delta_S[k]} implies that, for each $j\in \mathcal{N}$,
\begin{equation}\label{eq:calculation_of_deltasj}
    \Delta S_{j}=\max\left\{0,(r_{j}-\tilde{S}_{j}), \max_{j^{\prime}\in\mathcal{N}\setminus \{j\}}\{u(\tilde{\omega}_j-\tilde{\omega}_{j^{\prime}})(r_{j^{\prime}}-\tilde{S}_{j^{\prime}})\}\right\}.
\end{equation}

\begin{figure}[!t]
\begin{center}
        \includegraphics[width=0.8\linewidth]{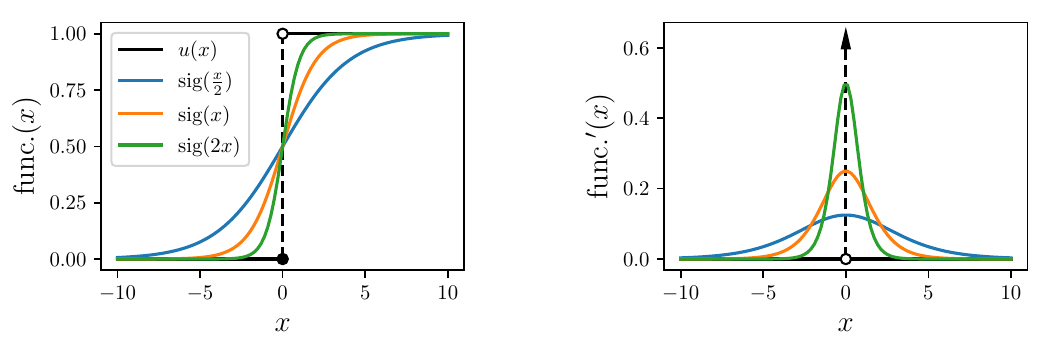}
    \end{center}
    \caption{Step function and sigmoid function}
    \label{Fig:sigmoid}
\end{figure}

\textbf{Smooth approximation:} In the above, the feasibility layer $F_{\mathcal{V}}(\cdot)$ is established by using the step function $u(\cdot)$ and the $\max$ operation. However, the gradient of the step function is zero at most places, as shown in the left part of Figure \ref{Fig:sigmoid}. We approximate the step function by  a smooth Sigmoid function as 
\begin{equation}
     u(\tilde{\omega}_{j}-\tilde{\omega}_{j^{\prime}}) \approx
    \mathrm{sig}(\varphi \cdot (\tilde{\omega}_{j}-\tilde{\omega}_{j^{\prime}})) = \frac{1}{1+\exp(-\varphi \cdot (\tilde{\omega}_{j}-\tilde{\omega}_{j^{\prime}}))}.
\end{equation}
The function is also shown in the right part of Figure \ref{Fig:sigmoid}.  In the above equation, $\varphi$ is a pre-specified  positive parameter. When $\varphi$ is a large enough value, the step function is well approximated by the sigmoid function, but the gradients at most places are almost zero. By setting $\varphi$ appropriately, the step function is properly approximated while having reasonable gradients. Similar to the dynamically adjusted target value in the step size \eqref{online_polyak}, we also dynamically adjust the value of $\varphi$. The details are given in the Appendix~\ref{app_pseudo_codes_online}.

The max operation within \eqref{eq:calculation_of_deltasj} makes $F_{\mathcal{V}}(\cdot)$ a piece-wise function. We do not try to smooth the $\max$ operation. At any non-differentiable point caused by the $\max$ operation, we simply use the right hand gradient of  the point as the gradient at the point. 

\section{Computational Results}\label{Sec:Numerical_Results}
In this section, we report the performance of our approach based on various sizes of test instances of various individual single-machine min-sum scheduling problems. The testing platform is equipped with AMD EPYC 7702, NVIDIA RTX 3090 GPU, Linux Ubuntu 18.04.5, and NVIDIA CUDA 12.4. Commercial solver IBM CPLEX 12.10 is used whenever there are LP or IP problems to be solved. The UMSNN, the supervised training, and the online single-instance learning are implemented by using Python 3.7 and Torch 1.13.1+cu117. Our code, datasets, and trained learnable parameters will be made available upon the publication of this study. 

In the following, we first describe in Section~\ref{sub:problem_type} which individual problems  and parameter distributions are tested. Then, we describe in Section~\ref{sub:generation_of_datasets} how the training, validation, and testing instances are generated. In Section~\ref{sec:benchmarks}, we describe two benchmark solution approaches to be compared with our approach. Finally, in Section~\ref{sec:results}, we report the computational results.

\subsection{Test Problems and Parameter Distributions} \label{sub:problem_type}
We test the performance of our approach based on the nine specific single-machine min-sum problems shown in Table~\ref{tab:instace_types}.
The test instances of these problems are generated following the distributions  of the problem parameters defined as follows:
\begin{itemize}
 \item Problem size, $n$, is drawn from the two size groups: size group 1, where $n\in \{500, 600, 700\}$, and size group 2, where $n\in \{800, 900, 1000\}$.   
 \item The parameters $p_j$, $w_j$, $w_j^A$, $w_j^B$, $w_{1j}$ and $w_{2j}$ are all drawn  randomly from U$\{1, 100\}$, where U$\{x, y\}$ denotes a discrete uniform distribution between $x$ and $y$.  For problems with release times, the release times are drawn from U$\{1,   0.5P\}$. For problems with due dates, the due dates $d_j$ are generated from U$\{1,   \xi^{d} P\}$, where parameter $\xi^{d}$ controls the tightness of the due dates, and $\xi^{d} \in \{0.2, 0.5, 0.8\}$ for some problems and $\xi^{d} = 0.5$ for some other problems. 
 \item The parameter $a_j$ in the objective function $\sum w_jC_j^{a_j}$ is drawn from U$[0.5, 1.5]$, where U$[x,y]$ is a continuous uniform distribution between $x$ and $y$. For the two bi-criterion problems, $\rho \in \{0.3, 0.7\}$. For the two two-agent problems, $\rho = 0.5$ and $|\mathcal{N}_A|/|\mathcal{N}_B|\in \{0.3, 0.7\}$. 
\end{itemize}
For each size group, Table~\ref{tab:instace_types} lists  the  17 problem cases and their corresponding parameter distributions to be tested. 
Consequently, there are 34  problem cases in total.

\begin{table}[ht]
    \centering 
    \caption{Problem cases and parameter distributions} 
    \label{tab:instace_types} 
    \begin{adjustbox}{width=\textwidth}
    
    \begin{tabular}{ccccccccc}  
        \toprule 
        Problem &Case& $p_j$&$d_j \sim U\{0,\xi^{d}P\}$&$r_j\sim U\{0,\xi^{r}P\}$& $w_j,w_{1j},w_{2j}$&Others\\
        \midrule 
        \multirow{3}{*}{$1||\sum w_j T_j $}&1&U\{1,100\}&$\xi^{d}=0.2$ &-&U\{1,100\}&-\\
          &2&U\{1,100\}&$\xi^{d}=0.5$ &-&U\{1,100\}&-\\
          &3&U\{1,100\}&$\xi^{d}=0.8$ &-&U\{1,100\}&-\\
        \midrule
        \multirow{2}{*}{$1||\sum \rho w_{1j}T_j+(1-\rho) w_{2j}C_j$}&4&U\{1,100\}&  $\xi^{d}=0.5$ &-&U\{1,100\}&$\rho=0.3$\\
          &5&U\{1,100\}&$\xi^{d}=0.5$ &-&U\{1,100\}&$\rho=0.7$\\
        \midrule
        \multirow{2}{*}{$1||\rho\sum w^{A}_jT^{A}_j+(1-\rho)\sum w^{B}_jC^{B}_j$}&6&U\{1,100\}&  $\xi^{d}=0.5$ &-&U\{1,100\}&$\rho=0.5, |\mathcal{N}_A|/|\mathcal{N}|=0.3$\\
        &7&U\{1,100\}&$\xi^{d}=0.5$ &-&U\{1,100\}&$\rho=0.5, |\mathcal{N}_A|/|\mathcal{N}|=0.7$\\
        \midrule
        $1||\sum w_j C_j^{a_j} $ &8&U\{1,100\}&-&-&U\{1,100\}&$a_j\sim U[0.5,1.5]$\\
        \midrule
        $1|r_j|\sum w_j C_j $ &9&U\{1,100\}&-&$\xi^{r}=0.5$ &U\{1,100\}&-\\
        \midrule
        \multirow{3}{*}{$1|r_j|\sum w_j T_j $}  &10&U\{1,100\}&$\xi^{d}=0.2$ &$\xi^{r}=0.5$ &U\{1,100\}&-\\
           &11&U\{1,100\}&$\xi^{d}=0.5$ &$\xi^{r}=0.5$ &U\{1,100\}&-\\
          &12&U\{1,100\}&$\xi^{d}=0.8$ &$\xi^{r}=0.5$ &U\{1,100\}&-\\
        \midrule
        \multirow{2}{*}{$1|r_j|\rho\sum w^{A}_jT^{A}_j+(1-\rho)\sum w^{B}_jC^{B}_j$}  &13&U\{1,100\}&$\xi^{d}=0.5$ &$\xi^{r}=0.5$ &U\{1,100\}&$\rho=0.5, |\mathcal{N}_A|/|\mathcal{N}|=0.3$\\
           &14&U\{1,100\}&$\xi^{d}=0.5$ &$\xi^{r}=0.5$ &U\{1,100\}&$\rho=0.5, |\mathcal{N}_A|/|\mathcal{N}|=0.7$\\
        \midrule
        \multirow{2}{*}{$1|r_j|\sum \rho w_{1j}T_j+(1-\rho) w_{2j}C_j$}&15&U\{1,100\}&  $\xi^{d}=0.5$ &$\xi^{r}=0.5$ &U\{1,100\}&$\rho=0.3$\\
           &16&U\{1,100\}&$\xi^{d}=0.5$ &$\xi^{r}=0.5$ &U\{1,100\}&$\rho=0.7$\\
        \midrule
        $1|r_j|\sum w_j C_j^{a_j} $  &17&U\{1,100\}&-&$\xi^{r}=0.5$ &U\{1,100\}&$a_j\sim U[0.5,1.5]$\\
        \bottomrule 
    \end{tabular}
    \end{adjustbox}
\end{table}

\subsection{Training, validation and testing sets}
\label{sub:generation_of_datasets}
As discussed in Section~\ref{sec:machine-learning}, in order for a supervised learning model to have a good generalization capability, the model needs to be trained  using input instances generated following the same or a similar distribution as the input space of the online instances to be solved. 
Therefore, we train our UMSNN model separately for each of the nine problems, except that for problems $1|r_j|\sum w_jC_j$ and $1|r_j|\sum w_jT_j$, the model is trained together. Moreover, the model is trained separately for each of the two size groups. However, the different parameter cases of a problem for each size group are trained together. For example,  the three cases of problem $1||\sum w_jT_j$ shown in Table~\ref{tab:instace_types} for problem sizes $n\in\{500,600,700\}$ are trained together, and these same three problem cases for problem sizes $n\in\{800,900,1000\}$ are also trained together, but separately from the smaller problem sizes. 

Ideally, the instances within a training set should follow exactly the same distribution (in terms of both the number of jobs and the parameters) as in  the instances to be solved online, which are described in the previous subsection. This, however, is impractical, because it would take long time to find optimal solutions  to training instances generated this way  due to their NP-hardness. Therefore, the training sets we use  consist of specially designed instances for which optimal solutions can be found quickly. The validation sets are used to periodically evaluate the model’s performance  to guide hyperparameter adjustments and determine when to stop training. The testing sets are used to evaluate the final performance of a trained model.  Both the validation sets and testing sets are generated following distributions described in the previous subsection. In the following, we describe how these data sets are generated precisely. 

\textbf{Generation of training sets:} 
 For each of the 17 problem cases shown in Table~\ref{tab:instace_types}, for each of the two size groups, a training set with 10k special instances   and the corresponding optimal solutions are generated as follows.  We first generate 2.5k small instances  with the number of jobs $n$ randomly drawn from U$\{500, 700\}$ or U$\{800, 1000\}$ depending on the  size group being considered,  the job processing times drawn from U$\{1,5\}$, and all the other parameters generated following the distributions given in Table~\ref{tab:instace_types} except that the weights $w_j$ for problems $1||\sum w_jC_j^{a_j}$ and $1|r_j|\sum w_jC_j^{a_j}$ are generated differently as described below. Optimal solutions to these instances are found by solving the time-indexed IP formulations of these instances using CPLEX. The small instances are scaled up to obtain large instances following the scheme shown in Table~\ref{tab:Approach to generating special instances}. Two values of the  scalar $\alpha$ are used:  $\alpha=20$, which makes the maximum job processing time equal to 100, and $\alpha$ is randomly selected from  U$\{1, 19\}$. This gives a total of 5k special instances, for which the corresponding optimal solutions are easily obtained based on the optimal solutions of the corresponding small instances. To enrich the training datasets, following the data augmentation approach presented in subsection \ref{subsec:data_augment}, we generate one extra instance corresponding to each  special instance generated above based on the optimal solution of the instance. For an instance  involving due dates, we randomly select half of the on-time jobs (i.e., $S_j+p_j < d_j$) in the optimal solution and change their due dates to new due dates drawn from U$\{S_j+p_j, d_j-1\}$. For an instance involving release times, we randomly select half of the jobs and change their release dates to new release times drawn from U$\{r_j+1, S_j\}$. After the data augmentation, a total of 10k special  instances and the corresponding optimal solutions are generated for each problem case.

For problems $1||\sum w_jC_j^{a_j}$ and $1|r_j|\sum w_jC_j^{a_j}$,  since following the scheme in Table~\ref{tab:Approach to generating special instances}, when a smaller instance is scaled up to a larger instance, the weights of the jobs need to be divided by $\alpha^{a_j}$, we generate the job weights from U$\{1,100\}\times \alpha^{a_j}$, which then ensures that the weights of the large instances to be generated fall within  U$\{1,100\}$.

\textbf{Generation of validation sets and testing sets:} 
For each problem case and each size group, a validation set of 15 instances is generated following the distributions given in  Table~\ref{tab:instace_types}. For the first size group, the 15 instances consist  of 5 with 500 jobs, 5 with 600 jobs and 5 with 700 jobs, and for the second size group, the 15 instances consist of 5 with 800 jobs, 5 with 900 jobs and 5 with 1000 jobs. Similarly, for each problem case and each size group, a testing set of 15 instances are generated the same way.

\subsection{Benchmarks} 
\label{sec:benchmarks}

The fact that the UMSNN is built on the time-indexed formulation \label{time-index-formulation} motivates us to utilize this same formulation to design the following two optimization based heuristics as benchmarks and use them to evaluate the performance of our approach. 

\textbf{Shrink + IP:}
The first benchmark is inspired by our idea in Section~\ref{subsec:special_instances} where we train the UMSNN using specially designed large instances, which are enlarged from randomly generated smaller instances by enlarging some parameters such as $p_j$, $r_j$, or / and $d_j$  following the scheme shown in Table \ref{tab:Approach to generating special instances}. The success of this scheme is partly due to the fact that the  time-indexed formulations of the smaller instances can be solved to optimality within a reasonable time. Now, following this idea, we reverse this process and ask the following question: for a randomly generated large instance, can we construct a smaller instance such that (i) the smaller instance can be solved to optimality quickly using the time-indexed formulation, and (ii) the optimal solution of the smaller instance can be easily expanded to become a feasible solution for the large instance? This can be done as follows. Given a large testing instance, we first generate a smaller instance by ``shrinking" the job processing times by 20 times as $\tilde{p}_j=\lceil p_j/20\rceil, \forall j\in \mathcal{N}$. The other parameters are shrunk accordingly following Table \ref{tab:Approach to generating special instances}. For example, for $1||\sum w_jT_j$, the due dates are calculated as $\tilde{d}_j=\lceil d_j/20\rceil, \forall j\in \mathcal{N}$. Then, we solve the time-indexed formulation of the smaller instance to obtain the best possible solution within a time limit following a warm-start strategy, as described below. Finally, a feasible solution to the given testing instance is generated by scheduling the jobs as tightly as possible using the same job sequence in the solution to the smaller instance. We call this benchmark approach ``Shrink + IP".

When using CPLEX to solve the time-indexed formulations of smaller instances, we try to reduce the required computational time by leveraging CPLEX’s warm-start functionality. Specifically, for a given smaller instance, we first solve the LP relaxation of its time-indexed formulation to optimality using CPLEX and then generate a feasible solution based on the sampling heuristic described in Section \ref{subsec:list heuristic}. This feasible solution is generally of high quality. We then provide this solution to CPLEX as a warm start for solving the original time-indexed formulation of the smaller instance. Numerical results suggest that CPLEX can quickly find a high-quality solution, although guaranteeing an optimal solution often requires significantly more time. To save computational effort, we limit the computation time   after providing a warm start to CPLEX to 600 seconds.

\textbf{LP + sampling:}
The second benchmark is a sampling heuristic similar to the one described in Section~\ref{subsec:list heuristic} except that we use the time-indexed formulation directly, instead of the UMSNN, to generate the probability vectors of the jobs' starting times. Given a randomly generated large testing instance, we first optimally solve the LP relaxation of the time-indexed formulation for this instance. In the solution obtained for job $j$,  $\{\tilde{x}_{jt}\}_{t\in\mathcal{T}}$  can be viewed as a probability vector of having time $t$ as the starting time of job $j$, since each element of it is between zero and one, and the sum of all the elements is 1. We then use the sampling heuristic based on the values of $\{\tilde{x}_{jt}\}_{j\in \mathcal{N}, t\in\mathcal{T}}$ given in Section \ref{subsec:list heuristic} to generate a feasible solution. We call this benchmark approach ``LP + sampling".

For the instances we tested on, which have very large sizes, it can be extremely time-consuming to solve the LP relaxation of the time-indexed formulation using the simplex method. We thus use the barrier crossover algorithm. 
We limit the computation time to 7200 seconds.

\subsection{Results}
\label{sec:results} 
We first describe how the supervised training is performed. Then we show the test results on the following three specific approaches, as compared to the two benchmark approaches described in Section~\ref{sec:benchmarks}: 
\begin{itemize}
\item ``Supervised + greedy", which takes the output from the trained model and applies the greedy approach described in Section~\ref{subsec:list heuristic} to find a feasible solution, 
\item ``Supervised + sampling", which is similar to ``supervised + greedy" except that the sampling approach described in Section~\ref{subsec:list heuristic} is used to find a feasible solution,
\item ``Supervised + online", which is the integrated offline and online learning approach, where  after the model is trained offline, online single-instance learning is applied.
\end{itemize} 

\textbf{Supervised training process:} The UMSNN model is trained following the same process described here for every problem. As described in Section \ref{sec:UMSNN}, our supervised learning objective is to solve problem \eqref{goal_of_offline} using the generated training sets. The layer sizes of the UMSNN are listed in Appendix~\ref{layer_sizes}. We adopt the widely used ADAM optimizer \citep{kingma2014adam}, a variant of stochastic gradient descent, with a batch size of 8. The learning rate starts at 0.00005 and decreases by a factor of 0.9 after each epoch, allowing for more aggressive early updates and more refined adjustments as training progresses. After every epoch, we evaluate the mean gap on the validation sets, and the greedy heuristic is used to generate a feasible solution for each prediction. If the mean gap fails to improve for two consecutive epochs, we halt training to prevent overfitting.

\begin{figure}[htbp]
    \centering
    \subfigure[Training Loss]{
        \includegraphics[width=0.41\textwidth]{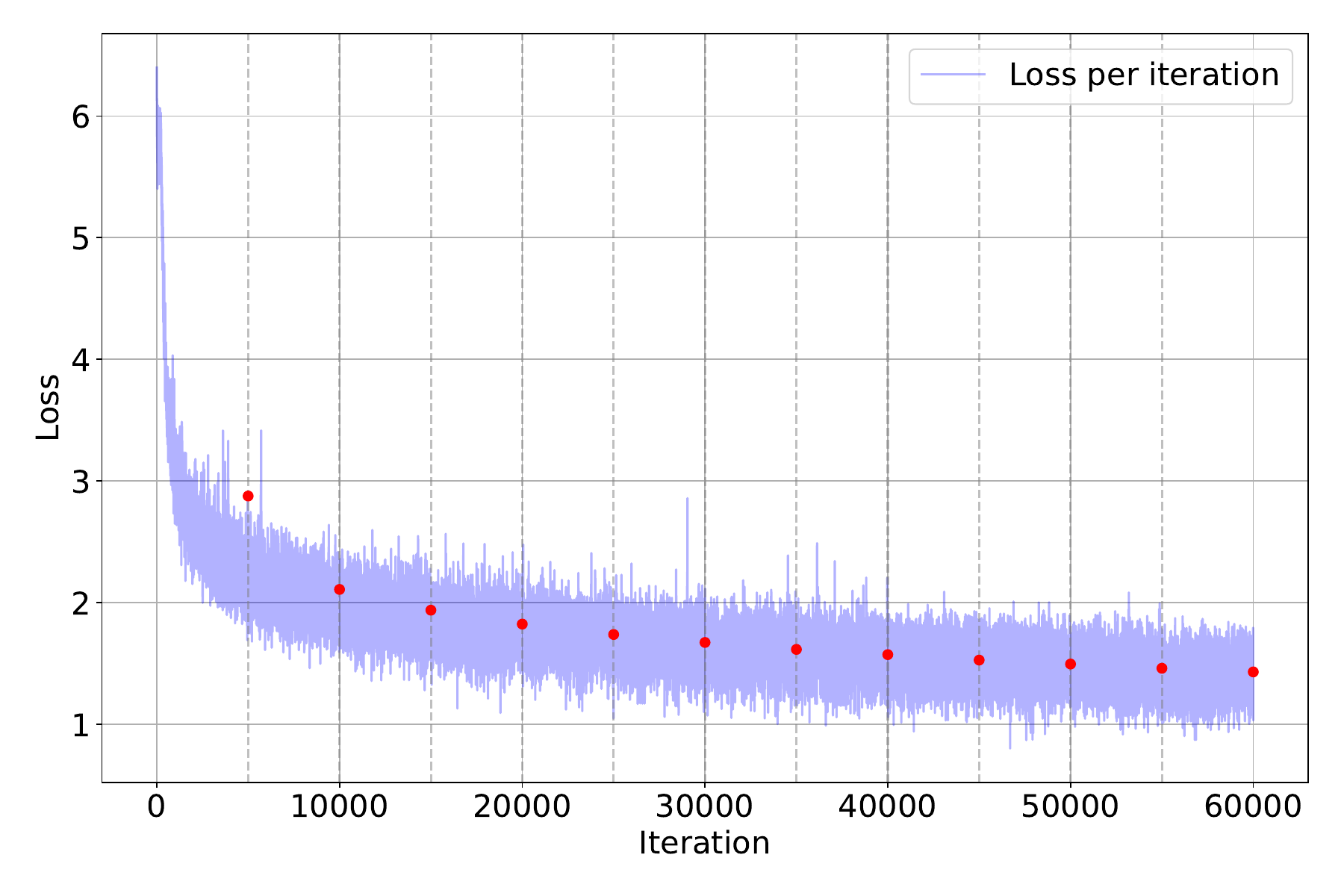} 
        \label{Fig_illustration_for_supervised_a}
    }
    \subfigure[Performance on the validation sets]{
        \includegraphics[width=0.45\textwidth]{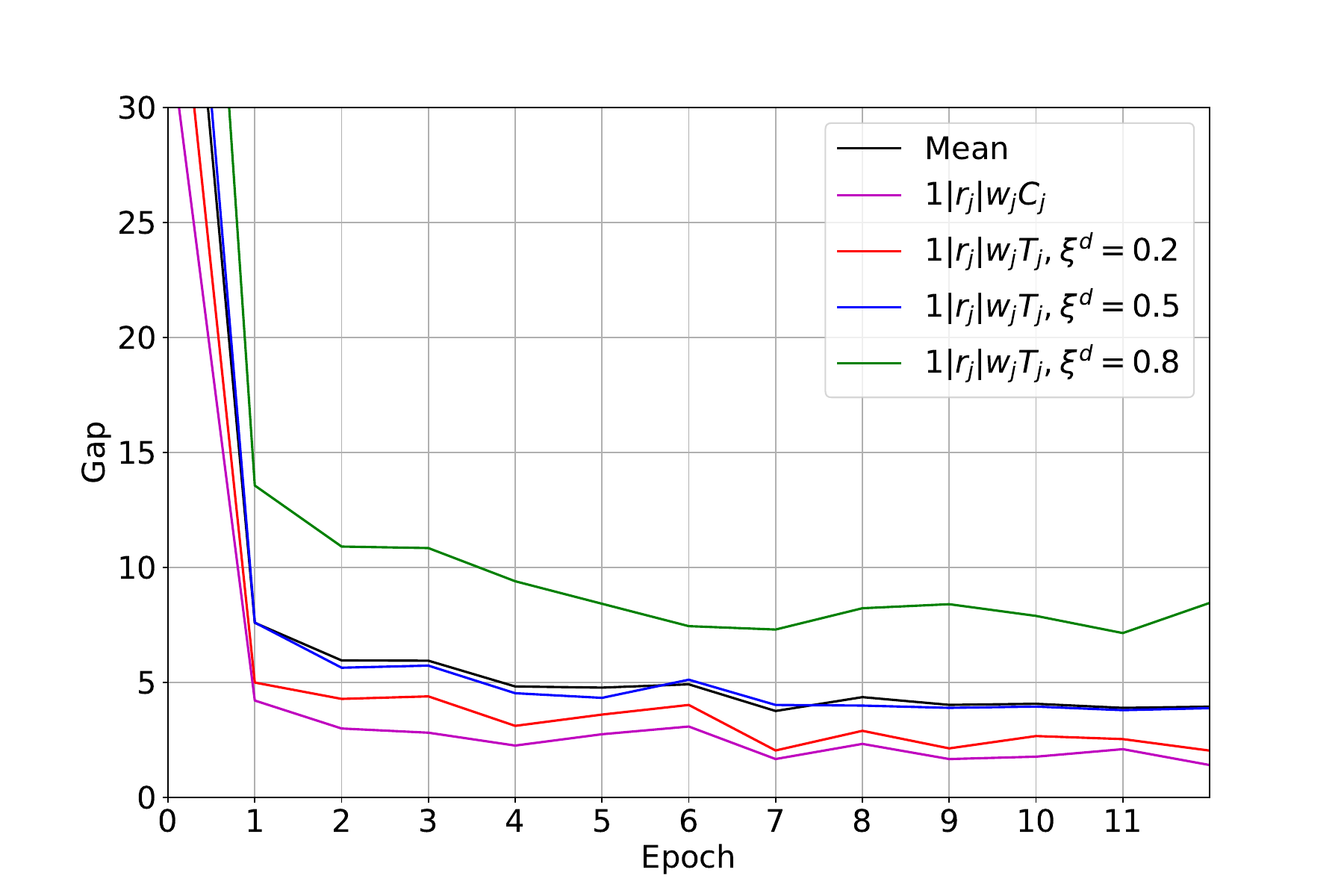}
        \label{Fig_illustration_for_supervised_b}
    }
    \caption{Offline supervised learning for subclasses M-9 to M-12}
    \label{Fig_illustration_for_supervised}
\end{figure}

To illustrate, we visualize the training process for problems $1|r_j|\sum w_jC_j$ and  $1|r_j|\sum w_jT_j$ with 500 to 700 jobs, by showing the losses and the validation accuracies. As discussed in Section~\ref{sub:generation_of_datasets}, all the associated four  cases of these problems (see Table~\ref{tab:instace_types}) are trained together. Since the  training set for each problem case contains 10,000  instances (as described in Section~\ref{sub:generation_of_datasets}), there are a total of 40,000 instances. Therefore, each training epoch involves 5,000 learnable parameter updates. Figure \ref{Fig_illustration_for_supervised_a} illustrates the  training process, showing the loss for each iteration (in blue) and the average loss per epoch (in red). The loss decreases notably during the first few epochs. The gaps on the validation sets are shown in Figure \ref{Fig_illustration_for_supervised_b}. Initially, before the learnable parameters are updated, the average gap exceeds 30\%. After the first epoch, it decreases sharply to about 7\%. The gap continues to decline over several subsequent epochs. By the 11th and 12th epochs, no improvement is observed, triggering the termination of the training process.
Thus, the learnable parameters obtained at the 12th epoch are used to solve the testing instances.

\textbf{Summary of the test results:} The test results of our three approaches ``Supervised + greedy", ``Supervised + sampling" and ``Supervised + online", and the two benchmark approaches ``LP + sampling" and ``Shrink + IP" for the two size groups are reported in Tables~\ref{tab:EX500TO700} and \ref{tab:EX800TO1000}, respectively. In these tables, each row involves 5 test instances.  
For each approach, the two columns under the measure ``Gap" show the average and maximum relative gap (in percentage) between the objective value of the solution found by this approach and the lower bound found by solving the LP relaxation of the time-indexed formulation. For every testing instance, approach ``Supervised + greedy"  always takes less than one second, and approach ``Supervised + sampling" always takes less than seven seconds to find a solution. Thus, computation times for these two approaches are not reported in the table. For approaches ``Supervised + online" and ``Shrink + IP", the two columns under ``Time" measure show the average and maximum computational time. Approach ``LP + sampling" can be very time-consuming, and hence for this approach, we limit the time that can be spent on each instance to 7200 seconds. For this approach, the column ``\#" reports the number of instances (out of 5) solved by this approach within this time limit, and the two columns under ``Time" measure show the average and maximum computational time for the instances solved within the time limit. 

Based on these results, we can make the following observations. 
\begin{itemize}
\item The two approaches based on the supervised learning only without online learning take very little time and hence are extremely efficient. They are also very effective because they generate  solutions with a less than 5\% gap from the lower bound for all instances except for the instances of the problems with the $\sum w_jT_j$ objective when the due dates are relatively loose. By comparing these two approaches, it is quite clear that the sampling heuristic almost always generates better solutions than the greedy heuristic. 
\item The integrated learning approach generates by far the best solutions among the three learning based approaches. In most cases, it generates a solution within 1 to 2\% of the lower bound. Compared to the other two learning approaches where online learning is not used,  the added online learning procedure in the integrated learning approach takes an extra 20 to 100 seconds for most instances. However, the online learning procedure improves the solution quality significantly.  Overall,  the integrated learning approach achieves a very good balance in terms of both solution effectiveness and computational efficiency. For an illustration of the online learning process, see Appendix~\ref{illustration_online}.
\item Although the benchmark approach ``LP + sampling" generates even better solutions than the integrated learning approach, it is very time-consuming, and for majority of the instances with 800 or more jobs, it cannot find a feasible solution in two hours. This reflects the strength of the time-indexed formulation (i.e., it has tight LP relaxation) as well as the weakness of such a formulation (i.e., it is large in scale and can take a long time to solve). 
\item The benchmark approach ``Shrink + IP" performs poorly in terms of both solution quality and computational efficiency. 
\end{itemize}

\begin{table}[h]
    \centering 
    \caption{Offline training results on generic instances with 500 to 700 jobs} 
    \label{tab:EX500TO700} 
    \begin{adjustbox}{width=\textwidth}
    \begin{tabular}{ccccccccccc}  
        \toprule 
        \multicolumn{2}{c}{Instances} & {\makecell{Supervised\\+greedy }}& {\makecell{Supervised\\+sampling }}& \multicolumn{2}{c}{Supervised+online} & \multicolumn{3}{c}{LP+Sampling}& \multicolumn{2}{c}{Shrink+IP}\\
        \cmidrule(lr){1-2}   \cmidrule(lr){5-6} \cmidrule(lr){7-9} \cmidrule(lr){10-11} 
        Type & $|\mathcal{N}|$ &  \makecell{Gap (\%)\\(avg - max)}& \makecell{Gap (\%)\\(avg - max)}& \makecell{Gap (\%)\\(avg - max)} &\makecell{Time (s)\\(avg - max)}&$\#$& \makecell{Gap (\%)\\(avg - max)} &\makecell{Time (s)\\(avg - max)}& \makecell{Gap (\%)\\(avg - max)} &\makecell{Time (s)\\(avg - max)}\\
        \midrule 
        
        \multirow{3}{*}{\makecell{$1||\sum w_j T_j$\\$ \xi^{d}=0.2$}}
        &500& 1.55 - 1.94 & 1.24 - 1.6  &0.19 - 0.4 & 23.5 - 32.4&5&0.04 - 0.06 & 1401 - 1477& 4.1 - 4.8 & 193.4 - 214.0\\
        &600& 1.34 - 1.43 & 0.96 - 1.18 &0.14 - 0.21 & 27.2 - 42.0 &5&0.03 - 0.04 & 1933 - 2618& 3.9 - 4.3 & 254.2 - 267.3\\
        &700& 1.17 - 1.32 & 0.92 - 1.04 &0.11 - 0.13 & 40.8 - 51.1&5 &0.04 - 0.04 & 2800 - 3570& 4.3 - 4.9 & 337.4 - 392.5\\
        \cmidrule(lr){1-11}
        \multirow{3}{*}{\makecell{$1||\sum w_j T_j$\\$ \xi^{d}=0.5$}} 
        &500& 1.89 - 3.04  & 1.25 - 2.06& 0.62 - 1.09 & 20.5 - 26.9&5&0.2 - 0.27 & 1228 - 1275 & 16.6 - 18.2 & 313.1 - 485.5\\
        &600& 1.89 - 2.42& 1.1 - 1.25&0.51 - 0.56 & 25.1 - 30.3&5&0.19 - 0.22 & 2948 - 6549& 16.8 - 19.7 & 569.7 - 810.7\\
        &700& 1.57 - 2.42& 1.02 - 1.16&0.39 - 0.5 & 38.0 - 44.9 &5&0.15 - 0.19 & 2729 - 2928& 15.4 - 15.9 & 698.2 - 842.5\\
        \cmidrule(lr){1-11}
        \multirow{3}{*}{\makecell{$1||\sum w_j T_j$\\$ \xi^{d}=0.8$}} 
        &500& 8.95 - 10.77& 5.47 - 6.67 &2.27 - 3.34 & 17.2 - 20.3&5&0.8 - 1.0 & 1100 - 1212& 49.5 - 53.9 & 615.6 - 758.4\\
        &600&  6.44 - 10.43 &4.41 - 6.13 &1.92 - 2.87 & 29.2 - 48.0&5&0.62 - 0.89 & 1875 - 2167& 48.4 - 55.0 & 751.6 - 800.9\\
        &700& 8.69 - 10.25 &5.85 - 6.97&1.76 - 2.07 & 27.4 - 39.7&5& 0.66 - 0.73 & 2873 - 3552& 50.6 - 51.1 & 846.8 - 878.0\\
        \midrule 
        \multirow{3}{*}{\makecell{$1||\rho\sum w^{A}_jT^{A}_j$\\$+(1-\rho)\sum w^{B}_jC^{B}_j$\\$|\mathcal{N}_A|/|\mathcal{N}|=0.3$ }}  
        &500&  0.48 - 0.63 &0.44 - 0.6&0.05 - 0.06 & 36.6 - 44.9&5&0.03 - 0.04 & 1176 - 1283& 3.7 - 3.9 & 201.7 - 215.7\\
        &600&  0.65 - 0.86 &0.45 - 0.68&0.07 - 0.18 & 41.9 - 57.8&5&0.03 - 0.04 & 1925 - 2429& 3.7 - 4.2 & 254.3 - 261.7\\
        &700&  1.13 - 1.48 &0.7 - 0.79&0.04 - 0.05 & 53.2 - 70.8&4&0.02 - 0.03 & 2852 - 3305& 3.6 - 4.0 & 334.7 - 364.4\\
        \cmidrule(lr){1-11}
        \multirow{3}{*}{\makecell{$1||\rho\sum w^{A}_jT^{A}_j$\\$+(1-\rho)\sum w^{B}_jC^{B}_j$\\$|\mathcal{N}_A|/|\mathcal{N}|=0.7$ }}  
        &500&  1.77 - 2.05 &1.49 - 1.74&0.24 - 0.3 & 31.5 - 33.3&5&0.11 - 0.14 & 1176 - 1215& 10.1 - 11.4 & 247.3 - 285.8\\
        &600&  1.75 - 2.57 &1.31 - 1.52&0.19 - 0.25 & 33.7 - 39.5&5&0.1 - 0.13 & 1739 - 1834& 9.7 - 10.9 & 297.9 - 349.8\\
        &700&  2.08 - 2.45 &1.69 - 1.98&0.18 - 0.23 & 40.9 - 45.1 &5&0.09 - 0.11 & 2900 - 3321& 9.4 - 10.5 & 381.5 - 481.0\\
        \midrule 
        \multirow{3}{*}{\makecell{$1||\sum \rho w_{1j}T_j+$\\$(1-\rho) w_{2j}C_j$\\$\rho=0.3$ }}  
        &500&  0.19 - 0.24 &0.17 - 0.23 &0.03 - 0.1 & 41.7 - 55.5&5&0.01 - 0.01 & 1321 - 1596&  1.9 - 2.3 & 399.4 - 402.1\\
        &600&  0.35 - 0.95 &0.17 - 0.22 &0.02 - 0.04 & 44.8 - 60.1&5&0.01 - 0.01 & 2390 - 3967& 2.1 - 2.5 & 484.9 - 494.0\\
        &700&  0.54 - 1.3  &0.3 - 0.4 &0.03 - 0.09 & 58.2 - 82.5&5&0.01 - 0.01 & 2758 - 3343& 2.0 - 2.2 & 599.9 - 612.9\\
        \cmidrule(lr){1-11}
         \multirow{3}{*}{\makecell{$1||\sum \rho w_{1j}T_j+$\\$(1-\rho) w_{2j}C_j$\\$\rho=0.7$ }} 
         &500&  1.7 - 1.99 &1.46 - 1.74 &0.11 - 0.15 & 34.9 - 47.0&5&0.04 - 0.05 & 1239 - 1329& 6.8 - 7.4 & 416.5 - 429.7\\
        &600&  1.44 - 1.75 &1.22 - 1.53&0.1 - 0.14 & 39.1 - 48.7&5&0.04 - 0.05 & 2372 - 3740& 6.6 - 6.9 & 522.2 - 545.7\\
        &700& 1.61 - 2.27 &1.29 - 1.64&0.09 - 0.13 & 47.8 - 50.4 &5&0.03 - 0.04 & 2837 - 3194& 6.2 - 6.7 & 638.5 - 646.8\\
        \midrule 
         \multirow{3}{*}{\makecell{$1||\sum w_{j}C_j^{a_j}$ }} 
         &500& 0.33 - 0.38&  0.62 - 0.74&$<$0.01&41.3 - 64.2 &5&$<$0.01&1344 - 1426& 2.2 - 2.5 & 700.6 - 727.1\\
        &600& 0.49 - 0.7 & 0.68 - 0.87  &$<$0.01&44.1 - 112.1&5& $<$0.01&1903 - 1953& 2.0 - 2.5 & 880.2 - 906.4\\
        &700& 0.72 - 1.28& 0.7 - 0.94   &$<$0.01&88.4 - 142.6&5&$<$0.01&3139 - 4197& 2.0 - 2.6 & 1086.9 - 1110.8\\
         \midrule
        \multirow{3}{*}{$1|r_j|\sum w_j C_j$}
        &500& 2.13 - 2.73  &1.15 - 1.44&1.0 - 1.21 & 12.2 - 17.1&5& 0.45 - 0.6 & 1669 - 1901& 11.7 - 13.2 & 400.9 - 526.8\\
        &600& 2.16 - 2.94  &1.24 - 1.4&1.07 - 1.25 & 14.5 - 18.8&5 &0.37 - 0.45 & 2580 - 2784& 11.5 - 12.1 & 775.1 - 810.0\\
        &700& 2.0 - 2.88  &1.07 - 1.37&0.72 - 0.99 & 30.1 - 41.9 &5&0.39 - 0.41 & 3683 - 4497& 11.8 - 12.6 & 854.1 - 880.3\\
         \midrule
         \multirow{3}{*}{\makecell{$1|r_j|\sum w_j T_j$\\$ \xi^{d}=0.2$}}
        &500&  2.48 - 2.69  &1.51 - 1.97&1.44 - 1.92 & 13.2 - 19.8&5&0.59 - 0.69 & 2018 - 3707& 14.1 - 15.2 & 418.2 - 632.1\\
        &600&  2.63 - 2.87 &1.5 - 1.82&1.17 - 1.83 & 22.8 - 40.6&5&0.53 - 0.67 & 2509 - 2655& 15.3 - 16.0 & 783.6 - 800.7\\
        &700& 2.51 - 3.14 &1.57 - 1.81&1.05 - 1.62 & 26.6 - 39.3&5&0.49 - 0.55 & 3646 - 4039& 14.8 - 15.9 & 852.5 - 877.8\\
        \cmidrule(lr){1-11}
         \multirow{3}{*}{\makecell{$1|r_j|\sum w_j T_j$\\$ \xi^{d}=0.5$}}
        &500&  4.08 - 4.82 &2.01 - 2.2&1.77 - 2.12 & 17.5 - 29.4 &5&0.78 - 0.83 & 1668 - 1748& 20.4 - 20.7 & 389.3 - 515.5\\
        &600&  3.91 - 4.56  &2.21 - 2.5 &1.66 - 2.03 & 21.5 - 28.3&5&0.87 - 1.03 & 3580 - 5131& 20.1 - 21.2 & 597.1 - 799.4\\
        &700& 3.4 - 4.55 &1.98 - 2.69&1.27 - 1.55 & 31.2 - 34.7&3&0.68 - 0.76 & 4136 - 5066& 20.2 - 22.2 & 843.3 - 865.6\\
        \cmidrule(lr){1-11}
         \multirow{3}{*}{\makecell{$1|r_j|\sum w_j T_j$\\$ \xi^{d}=0.8$}}
        &500&  7.23 - 9.32  &4.39 - 5.2 &3.13 - 3.79 & 18.5 - 32.3&5&1.51 - 1.81 & 1558 - 1763&  24.9 - 26.7 & 615.5 - 758.0\\
        &600&  8.06 - 11.46 &3.98 - 4.69 &3.03 - 3.53 & 26.5 - 36.2&5&1.79 - 2.02 & 2453 - 2719& 25.8 - 27.6 & 765.0 - 813.6\\
        &700&  6.15 - 7.04  &3.2 - 3.73&2.31 - 2.78 & 28.9 - 39.7&5&1.47 - 1.66 & 3706 - 4093& 24.3 - 25.2 & 858.2 - 870.0\\
        \midrule
        \multirow{3}{*}{\makecell{$1|r_j|\rho\sum w^{A}_jT^{A}_j$\\$+(1-\rho)\sum w^{B}_jC^{B}_j$\\$|\mathcal{N}_A|/|\mathcal{N}|=0.3$ }} 
        &500& 2.64 - 3.64 &1.58 - 2.15&1.34 - 2.15 & 29.8 - 73.0& 5&0.48 - 0.56 & 1644 - 1811& 12.9 - 13.4 & 749.3 - 762.4\\
        &600& 2.81 - 3.24 &1.67 - 1.95 &1.52 - 1.87 & 16.8 - 26.4&5&0.48 - 0.54 & 2489 - 2721& 13.5 - 14.3 & 753.6 - 797.7\\
        &700& 1.81 - 2.33 &1.2 - 1.4&0.93 - 1.04 & 25.9 - 34.9 &5&0.45 - 0.59 & 3950 - 4652& 13.6 - 13.9 & 847.5 - 851.0\\
        \cmidrule(lr){1-11}
        \multirow{3}{*}{\makecell{$1|r_j|\rho\sum w^{A}_jT^{A}_j$\\$+(1-\rho)\sum w^{B}_jC^{B}_j$\\$|\mathcal{N}_A|/|\mathcal{N}|=0.7$ }} 
        &500&  3.91 - 4.48& 2.39 - 2.84 &1.89 - 2.45 & 17.6 - 28.2&5&0.68 - 0.81 & 1622 - 1888& 15.5 - 16.4 & 503.0 - 576.1\\
        &600&  2.8 - 4.09 &1.86 - 2.25&1.43 - 1.77 & 23.6 - 40.3&5&0.61 - 0.75 & 2577 - 2658& 16.0 - 16.6 & 798.4 - 799.9\\
        &700& 3.63 - 4.39 & 2.2 - 3.22&1.9 - 2.43 & 20.0 - 28.5 &3&0.62 - 0.72 & 3908 - 4446& 16.7 - 18.2 & 733.7 - 848.3\\
        \midrule
        \multirow{3}{*}{\makecell{$1|r_j|\sum \rho w_{1j}T_j+$\\$(1-\rho) w_{2j}C_j$\\$\rho=0.3$ }} 
        &500&  2.15 - 2.6  &1.23 - 1.59&0.88 - 1.09 & 20.2 - 32.7&5&0.4 - 0.47 & 1777 - 1985& 14.1 - 15.0 & 643.1 - 719.6\\
        &600&   1.79 - 2.2  &1.18 - 1.53&0.79 - 0.93 & 23.4 - 31.3&5&0.37 - 0.45 & 3019 - 3892& 14.1 - 14.6 & 1015.8 - 1050.7\\
        &700& 2.02 - 3.17  &1.47 - 2.0&1.07 - 1.5 & 23.1 - 25.9&5&0.38 - 0.4 & 4047 - 4289& 14.0 - 14.8 & 1115.4 - 1149.7\\
        \cmidrule(lr){1-11}
        \multirow{3}{*}{\makecell{$1|r_j|\sum \rho w_{1j}T_j+$\\$(1-\rho) w_{2j}C_j$\\$\rho=0.7$ }} 
        &500&  2.84 - 3.41& 1.82 - 2.11 &1.42 - 1.99 & 15.4 - 22.7&5&0.42 - 0.59 & 1937 - 2175& 16.9 - 18.1 & 806.2 - 942.3\\
        &600&  3.12 - 4.01& 1.78 - 2.55 &1.08 - 1.44 & 21.9 - 34.8&5&0.47 - 0.53 & 2783 - 3581& 17.5 - 18.3 & 958.2 - 1053.0\\
        &700& 2.62 - 3.24& 1.8 - 2.27 &1.13 - 1.34 & 28.5 - 34.6&4&0.46 - 0.53 & 4918 - 6898& 17.2 - 19.6 & 1132.1 - 1145.7\\
        \midrule
        \multirow{3}{*}{\makecell{$1|r_j|\sum w_{j}C_j^{a_j}$ }}  
        &500&   2.28 - 3.47 & 1.23 - 1.44   &0.85 - 1.01 & 26.7 - 36.6&5&0.55 - 0.65 & 2083 - 2177& 1.1 - 1.2 & 697.1 - 958.2\\
        &600&  1.72 - 2.3 & 1.09 - 1.31     &0.85 - 1.02 & 18.8 - 25.3&5&0.5 - 0.54 & 4124 - 6418& 1.1 - 1.3 & 1306.8 - 1446.9\\
        &700&  1.89 - 2.81 & 1.14 - 1.4     &0.95 - 1.17 & 23.5 - 49.2&5&0.57 - 0.67 & 4926 - 6689& 1.0 - 1.2 & 1373.1 - 1586.0\\
        \bottomrule 
    \end{tabular}
    \end{adjustbox}
\end{table}

\begin{table}[t]
    \centering 
    \caption{Offline training results on generic instances with 800 to 1000 jobs} 
    \label{tab:EX800TO1000} 

    \begin{adjustbox}{width=\textwidth}
    \begin{tabular}{ccccccccccc}  
        \toprule 
        \multicolumn{2}{c}{Instances} & {\makecell{Supervised\\+greedy }}& {\makecell{Supervised\\+sampling }}& \multicolumn{2}{c}{Supervised+online} & \multicolumn{3}{c}{LP+sampling }& \multicolumn{2}{c}{Shrink+IP}\\
        \cmidrule(lr){1-2}   \cmidrule(lr){5-6} \cmidrule(lr){7-9} \cmidrule(lr){10-11}
        Type & $|\mathcal{N}|$ &  \makecell{\makecell{Gap (\%)\\(avg - max)}}& \makecell{\makecell{Gap (\%)\\(avg - max)} }& \makecell{Gap (\%)\\(avg - max)} &\makecell{Time (s)\\(avg - max)}& $\#$&\makecell{Gap (\%)\\( avg - max)}&\makecell{Time (s)\\(avg - max)}&\makecell{Gap (\%)\\( avg - max)}&\makecell{Time (s)\\(avg - max)}\\
        \midrule 
         \multirow{3}{*}{\makecell{$1||\sum w_j T_j$\\$ \xi^{d}=0.2$}}
        &800&   1.36 - 1.51 & 0.94 - 1.1 &0.12 - 0.22 & 37.3 - 62.8   &5&0.03 - 0.04 & 4166 - 4708& 4.1 - 4.6 & 421.8 - 521.2 \\
        &900&   1.07 - 1.58 & 0.77 - 1.05 &0.11 - 0.22 & 53.4 - 79.8 &3&0.02 - 0.03 & 4654 - 4829 & 4.4 - 4.8 & 494.4 - 524.3\\
        &1000&  1.09 - 1.31 & 0.83 - 1.02  & 0.06 - 0.09 & 73.1 - 93.6&4& 0.02 - 0.03 & 6310 - 6748& 4.3 - 4.8 & 576.8 - 601.3 \\
        \cmidrule(lr){1-11} 
         \multirow{3}{*}{\makecell{$1||\sum w_j T_j$\\$ \xi^{d}=0.5$}}
        &800&    2.84 - 3.56 & 1.56 - 1.92&0.4 - 0.49 & 38.4 - 66.2  & 5& 0.14 - 0.18 & 3988 - 4259 & 16.3 - 18.2 & 877.2 - 898.9\\
        &900&   2.75 - 4.34& 1.64 - 2.44 &0.31 - 0.39 & 56.1 - 83.2 &3 &0.15 - 0.17 & 5629 - 7029& 17.1 - 18.1 & 942.0 - 945.0 \\
        &1000&   2.8 - 3.73& 1.58 - 1.79 &0.37 - 0.44 & 51.7 - 62.1  &2 &0.12 - 0.13 & 7073 - 7195& 16.7 - 17.4 & 995.4 - 1008.5 \\
        \cmidrule(lr){1-11} 
         \multirow{3}{*}{\makecell{$1||\sum w_j T_j$\\$ \xi^{d}=0.8$}}
        &800&   6.01 - 10.86& 3.32 - 3.76&1.35 - 1.46 & 32.0 - 38.6 &3&0.54 - 0.57 & 3601 - 4036 & 52.7 - 56.3 & 895.4 - 901.6\\
        &900&   7.46 - 8.89& 4.16 - 4.83&1.43 - 1.84 & 39.2 - 57.7    &3&0.56 - 0.62 & 4817 - 5056 & 54.0 - 56.8 & 945.6 - 955.5\\
        &1000&   9.31 - 12.96& 4.83 - 5.72&1.54 - 1.83 & 57.4 - 87.9&1&0.66 - 0.66 & 6057 - 6057 & 54.5 - 58.3 & 998.6 - 1004.5\\
        \midrule
        \multirow{3}{*}{\makecell{$1||\rho\sum w^{A}_jT^{A}_j$\\$+(1-\rho)\sum w^{B}_jC^{B}_j$\\$|\mathcal{N}_A|/|\mathcal{N}|=0.3$ }} 
        &800&  1.05 - 1.42& 0.97 - 1.24 &0.08 - 0.17 & 49.4 - 80.1 &5&0.02 - 0.03 & 4100 - 4675& 4.2 - 4.8 & 427.4 - 450.9\\
        &900&  0.81 - 0.87& 0.74 - 0.78 &0.07 - 0.14 & 56.3 - 80.4 &3&0.02 - 0.02 & 6073 - 6592& 3.8 - 3.9 & 461.7 - 506.9\\
        &1000& 1.34 - 1.87& 1.12 - 1.29 &0.04 - 0.05 & 74.9 - 91.1 &3&0.02 - 0.02 & 6985 - 7115& 3.9 - 4.2 & 635.5 - 665.2\\
        \cmidrule(lr){1-11} 
        \multirow{3}{*}{\makecell{$1||\rho\sum w^{A}_jT^{A}_j$\\$+(1-\rho)\sum w^{B}_jC^{B}_j$\\$|\mathcal{N}_A|/|\mathcal{N}|=0.7$ }} 
        &800&  3.23 - 3.70& 3.04 - 3.34 &0.26 - 0.29 & 42.9 - 53.9&5&0.09 - 0.1 & 4152 - 4532& 9.7 - 11.2 & 646.3 - 821.1\\
        &900&  2.73 - 3.29& 2.41 - 2.85 &0.24 - 0.34 & 52.9 - 73.5 &4&0.07 - 0.08 & 6214 - 6755& 9.3 - 10.0 & 737.3 - 944.9\\
        &1000& 3.13 - 3.96& 2.56 - 3.33 &0.2 - 0.29 & 60.6 - 80.7 &1&0.06 - 0.06 & 6882 - 6882& 8.7 - 9.5 & 995.0 - 1002.8\\
        \midrule
        \multirow{3}{*}{\makecell{$1||\sum \rho w_{1j}T_j+$\\$(1-\rho) w_{2j}C_j$\\$\rho=0.3$ }} 
        &800&  0.31 - 0.39& 0.29 - 0.34 &0.01 - 0.01 & 108.5 - 127.1 & 2&0.01 - 0.01 & 5038 - 6280 & 2.0 - 2.2 & 851.6 - 908.5\\
        &900& 0.4 - 0.65&0.27 - 0.32 &0.01 - 0.01 & 98.5 - 110.1&1&0.01 - 0.01 & 5760 - 5760 & 1.9 - 2.1 & 1207.1 - 1952.8\\
        &1000& 0.87 - 1.19&0.49 - 0.61&0.01 - 0.01 & 112.2 - 131.5&0&-&- & 2.0 - 2.1 & 1145.7 - 1219.1\\
        \cmidrule(lr){1-11} 
        \multirow{3}{*}{\makecell{$1||\sum \rho w_{1j}T_j+$\\$(1-\rho) w_{2j}C_j$\\$\rho=0.7$ }} 
        &800&  2.04 - 2.39& 1.82 - 2.19 &0.09 - 0.11 & 80.2 - 90.6&3& 0.03 - 0.03 & 4560 - 6158& 6.9 - 7.6 & 925.9 - 1074.0 \\
        &900& 1.58 - 1.8 & 1.34 - 1.52 &0.07 - 0.09 & 76.0 - 102.6&5&0.03 - 0.03 & 5997 - 6674& 6.7 - 6.9 & 1281.4 - 2089.0 \\
        &1000&1.06 - 1.28& 0.89 - 1.23 &0.08 - 0.1 & 87.2 - 93.7 & 1&0.03 - 0.03 & 6958 - 6958 & 6.7 - 7.2 & 1252.9 - 1366.5\\
        \midrule
        \multirow{3}{*}{\makecell{$1||\sum w_{j}C_j^{a_j}$ }}  
        &800& 0.33 - 0.4& 0.4 - 0.46 &$<$0.01&178.9 - 226.9&4&$<$0.01&3766 - 4435& 2.0 - 2.3 & 1261.8 - 1289.9\\
        &900&  0.7 - 0.92& 0.6 - 0.84& $<$0.01&191.3 - 297.7&1&$<$0.01&5438 - 5438& 1.8 - 2.0 & 1456.9 - 1516.4\\
        &1000& 1.44 - 1.63& 0.75 - 0.95& $<$0.01& 243.3 - 342.0&0&-&- & 2.1 - 2.6 & 1486.2 - 1633.1\\
        \midrule
        \multirow{3}{*}{\makecell{$1|r_j|\sum w_j C_j$}}
        &800&    1.7 - 2.18 & 1.04 - 1.66 & 0.68 - 0.95 & 34.0 - 49.71 &4&0.33 - 0.37 & 5868 - 6206& 11.8 - 12.9 & 896.5 - 903.4\\
        &900&   1.5 - 2.19  & 0.82 - 0.97   &  0.79 - 0.96 & 31.1 - 56.45 &1&0.38 - 0.38 & 6826 - 6826& 11.8 - 11.9 & 945.6 - 957.2\\
        &1000&  1.62 - 2.03 & 1.05 - 1.39 &  0.82 - 1.15 & 34.3 - 54.85 &0&-&-& 11.6 - 12.0 & 993.6 - 1000.8\\
        
        \midrule 
         \multirow{3}{*}{\makecell{$1|r_j|\sum w_j T_j$\\$ \xi^{d}=0.2$}}
        &800&   2.45 - 3.0 &1.67 - 2.17 & 1.03 - 1.25 & 40.1 - 58.92&4&0.44 - 0.52 & 5801 - 6555& 15.6 - 16.1 & 898.2 - 901.7\\
        &900&   2.27 - 2.68  &1.4 - 1.91&1.05 - 1.31 & 34.7 - 50.45&0&-&-& 15.3 - 16.1 & 940.1 - 953.4\\
        &1000&  2.01 - 2.99 &1.27 - 1.74&0.83 - 1.3 & 52.0 - 80.77&0&-&-& 15.2 - 16.3 & 990.9 - 995.5\\
        \cmidrule(lr){1-11} 
         \multirow{3}{*}{\makecell{$1|r_j|\sum w_j T_j$\\$ \xi^{d}=0.5$}}
        &800&   3.56 - 4.47& 2.27 - 3.28&1.89 - 2.13 & 24.6 - 30.68&2&0.84 - 0.9 & 5923 - 6301& 19.7 - 20.7 & 895.2 - 904.1\\ 
        &900&   3.76 - 5.0 & 2.37 - 2.83&1.19 - 1.33 & 49.5 - 66.4&0&-&-& 20.2 - 20.5 & 946.4 - 955.6\\
        &1000&  3.5 - 5.0 & 2.35 - 3.44 &1.53 - 1.95 & 42.6 - 62.52&0&-&-& 19.9 - 20.7 & 997.4 - 1001.8\\
        \cmidrule(lr){1-11} 
         \multirow{3}{*}{\makecell{$1|r_j|\sum w_j T_j$\\$ \xi^{d}=0.8$}}
        &800&  7.34 - 9.91 &4.42 - 5.45&3.4 - 5.11 & 30.7 - 47.14 &4&1.52 - 1.67 & 5183 - 5605& 25.4 - 28.1 & 896.0 - 902.1\\
        &900&   5.87 - 8.7  &4.01 - 5.4 &2.98 - 4.39 & 45.9 - 72.75&0&-&-& 24.2 - 25.1 & 950.5 - 969.0\\
        &1000&    6.57 - 9.8 & 4.1 - 4.34 &2.61 - 4.0 & 65.3 - 84.4&0&-&-& 25.3 - 25.7 & 961.2 - 1032.1 \\
        \midrule
        \multirow{3}{*}{\makecell{$1|r_j|\rho\sum w^{A}_jT^{A}_j$\\$+(1-\rho)\sum w^{B}_jC^{B}_j$\\$|\mathcal{N}_A|/|\mathcal{N}|=0.3$ }} 
        &800&   2.38 - 3.27 &1.54 - 2.37 & 0.98 - 1.28 & 27.8 - 38.79&3&0.42 - 0.44 & 5637 - 6140& 13.4 - 14.0 & 894.2 - 907.2\\
        &900&   1.52 - 1.91 &1.04 - 1.15&0.94 - 1.22 & 37.0 - 85.78&0&-&-& 13.1 - 13.9 & 941.0 - 955.6\\
        &1000&   1.82 - 2.19 & 1.13 - 1.35&0.99 - 1.35 & 42.0 - 90.82&0&-&-& 13.4 - 14.3 & 998.0 - 1002.8\\
        \cmidrule(lr){1-11} 
        \multirow{3}{*}{\makecell{$1|r_j|\rho\sum w^{A}_jT^{A}_j$\\$+(1-\rho)\sum w^{B}_jC^{B}_j$\\$|\mathcal{N}_A|/|\mathcal{N}|=0.7$ }} 
        &800&  2.5 - 2.61 &  1.58 - 1.97 &1.25 - 1.58 & 33.6 - 56.77&4&0.53 - 0.58 & 6446 - 7162& 16.2 - 17.3 & 896.2 - 904.5\\
        &900&   2.77 - 3.2 &   1.85 - 2.34  &1.25 - 1.35 & 63.6 - 85.51&1&0.63 - 0.63 & 6717 - 6717& 16.0 - 16.7 & 946.5 - 972.7\\
        &1000&   3.55 - 4.64 & 1.72 - 2.14 &1.28 - 1.55 & 57.9 - 84.0&0&-&-& 16.2 - 17.4 & 999.0 - 1001.2\\
        \midrule
        \multirow{3}{*}{\makecell{$1|r_j|\sum \rho w_{1j}T_j+$\\$(1-\rho) w_{2j}C_j$\\$\rho=0.3$ }} 
        &800&   2.11 - 2.87&1.24 - 1.62 &0.9 - 1.24 & 33.7 - 56.07 &5&0.34 - 0.39 & 5834 - 6670 & 14.1 - 14.9 & 1532.2 - 2321.8\\
        &900&   1.77 - 2.15&1.27 - 1.46 &0.8 - 1.1 & 36.4 - 55.04&0&-&-  & 14.2 - 14.6 & 1420.0 - 1457.7\\
        &1000&   2.43 - 2.76 &1.72 - 2.03&1.04 - 1.49 & 40.8 - 53.27&0&-&-  & 14.0 - 14.3 & 1590.9 - 1664.4\\
        \cmidrule(lr){1-11} 
        \multirow{3}{*}{\makecell{$1|r_j|\sum \rho w_{1j}T_j+$\\$(1-\rho) w_{2j}C_j$\\$\rho=0.7$ }} 
        &800&  2.8 - 4.5 & 2.24 - 3.07 &1.31 - 2.0 & 28.1 - 35.5  &5&0.36 - 0.46 & 6297 - 6948 & 16.7 - 17.8 & 1537.9 - 2329.2\\
        &900&   3.58 - 4.38 & 2.55 - 3.21 &1.32 - 1.82 & 40.9 - 67.23 &0&-&-& 17.0 - 17.6 & 1419.0 - 1455.0 \\
        &1000&   2.9 - 4.43  &2.26 - 3.15 &1.29 - 1.8 & 46.1 - 87.18&0&-&- & 16.4 - 17.7 & 1591.2 - 1690.0\\
        \midrule
        \multirow{3}{*}{\makecell{$1|r_j|\sum w_{j}C_j^{a_j}$ }}  
        &800&    1.85 - 2.45 &1.08 - 1.19  & 0.85 - 1.03 & 37.9 - 46.19 &3&0.52 - 0.57 & 6419 - 6900 & 1.1 - 1.4 & 1713.4 - 1785.6\\
        &900&    2.42 - 3.11& 1.18 - 1.3 & 0.83 - 1.03 & 36.2 - 48.04 &0&-&-& 1.1 - 1.3 & 1740.6 - 1916.1\\
        &1000&    1.62 - 2.05&  1.03 - 1.19  &0.75 - 0.87 & 56.0 - 66.97 &0&-&-& 1.0 - 1.3 & 1822.2 - 1832.4\\
        \bottomrule 
    \end{tabular}

    \end{adjustbox}
\end{table}

\section{Conclusion}
\label{sec:conclusion}
In this paper, by exploiting the time-indexed formulation, we have proposed a unified supervised learning based approach  for the entire class of the single-machine problems with a min-sum objective. Through offline supervised training by using special instances, high-quality solutions are generated for instances with up to 1000 jobs in seconds. The added online single-instance learning procedure improves solution qualities significantly with an extra computational time of less than 100 seconds in most cases.  These advancements have major implications for scheduling in manufacturing, supply chain management, and beyond. Looking ahead, developing integrated offline and online training strategies for more complex scheduling problems, such as job shop scheduling, and for vehicle routing problems, represents a promising direction for future research.

\singlespacing
\setlength\bibsep{0pt}
\bibliographystyle{informs2014}
\bibliography{reference}

\begin{thebibliography}{75}
\providecommand{\natexlab}[1]{#1}
\providecommand{\url}[1]{\texttt{#1}}
\providecommand{\urlprefix}{URL }

\bibitem[{Agnetis et~al.(2014)Agnetis, Billaut, Gawiejnowicz, Pacciarelli,
  Soukhal et~al.}]{agnetis2014multiagent}
Agnetis A, Billaut JC, Gawiejnowicz S, Pacciarelli D, Soukhal A, et~al. (2014)
  Multiagent scheduling. \emph{Berlin Heidelberg: Springer Berlin Heidelberg.
  doi} 10(1007):978--3.

\bibitem[{Akturk \protect\BIBand{} Ozdemir(2001)}]{akturk2001new}
Akturk MS, Ozdemir D (2001) A new dominance rule to minimize total weighted
  tardiness with unequal release dates. \emph{European Journal of Operational
  Research} 135(2):394--412.

\bibitem[{Amari(1993)}]{amari1993backpropagation}
Amari Si (1993) Backpropagation and stochastic gradient descent method.
  \emph{Neurocomputing} 5(4-5):185--196.

\bibitem[{Anderson et~al.(1997)Anderson, Glass, Potts
  et~al.}]{anderson1997machine}
Anderson EJ, Glass CA, Potts CN, et~al. (1997) Machine scheduling. \emph{Local
  search in combinatorial optimization} 11:361--414.

\bibitem[{Avella et~al.(2005)Avella, Boccia, \protect\BIBand{}
  D’Auria}]{avella2005near}
Avella P, Boccia M, D’Auria B (2005) Near-optimal solutions of large-scale
  single-machine scheduling problems. \emph{INFORMS Journal on Computing}
  17(2):183--191.

\bibitem[{Ba et~al.(2016)Ba, Kiros, \protect\BIBand{} Hinton}]{ba2016layer}
Ba JL, Kiros JR, Hinton GE (2016) Layer normalization. \emph{arXiv preprint
  arXiv:1607.06450} .

\bibitem[{Balcan et~al.(2018)Balcan, Dick, Sandholm, \protect\BIBand{}
  Vitercik}]{balcan2018learning}
Balcan MF, Dick T, Sandholm T, Vitercik E (2018) Learning to branch.
  \emph{International conference on machine learning}, 344--353 (PMLR).

\bibitem[{Bengio et~al.(2021)Bengio, Lodi, \protect\BIBand{}
  Prouvost}]{bengio2021machine}
Bengio Y, Lodi A, Prouvost A (2021) Machine learning for combinatorial
  optimization: a methodological tour d’horizon. \emph{European Journal of
  Operational Research} 290(2):405--421.

\bibitem[{Berghman \protect\BIBand{} Spieksma(2015)}]{berghman-2015}
Berghman L, Spieksma FC (2015) Valid inequalities for a time-indexed
  formulation. \emph{Operations Research Letters} 43(3):268--272.

\bibitem[{Bilge et~al.(2007)Bilge, Kurtulan, \protect\BIBand{}
  K{\i}ra{\c{c}}}]{bilge2007tabu}
Bilge {\"U}, Kurtulan M, K{\i}ra{\c{c}} F (2007) A tabu search algorithm for
  the single machine total weighted tardiness problem. \emph{European Journal
  of Operational Research} 176(3):1423--1435.

\bibitem[{Blazewicz et~al.(2013)Blazewicz, Ecker, Pesch, Schmidt,
  \protect\BIBand{} Weglarz}]{blazewics-etal-2013}
Blazewicz J, Ecker KH, Pesch E, Schmidt G, Weglarz J (2013) \emph{Scheduling
  Computer and Manufacturing Processes} (springer science \& Business media).

\bibitem[{Bou{{\v{s}}}ka et~al.(2023)Bou{{\v{s}}}ka, {{\v{S}}}{{\u{u}}}cha,
  Nov{{\'a}}k, \protect\BIBand{} Hanz{{\'a}}lek}]{bouvska2023deep}
Bou{{\v{s}}}ka M, {{\v{S}}}{{\u{u}}}cha P, Nov{{\'a}}k A, Hanz{{\'a}}lek Z
  (2023) Deep learning-driven scheduling algorithm for a single machine problem
  minimizing the total tardiness. \emph{European Journal of Operational
  Research} 308(3):990--1006.

\bibitem[{Burke \protect\BIBand{} Curtois(2014)}]{burke2014new}
Burke EK, Curtois T (2014) New approaches to nurse rostering benchmark
  instances. \emph{European Journal of Operational Research} 237(1):71--81.

\bibitem[{Chang et~al.(2006)Chang, Hsieh, \protect\BIBand{}
  Liu}]{chang2006case}
Chang PC, Hsieh JC, Liu CH (2006) A case-injected genetic algorithm for single
  machine scheduling problems with release time. \emph{International Journal of
  Production Economics} 103(2):551--564.

\bibitem[{Cheung et~al.(2017)Cheung, Mestre, Shmoys, \protect\BIBand{}
  Verschae}]{cheung2017primal}
Cheung M, Mestre J, Shmoys DB, Verschae J (2017) A primal-dual approximation
  algorithm for min-sum single-machine scheduling problems. \emph{SIAM Journal
  on Discrete Mathematics} 31(2):825--838.

\bibitem[{Chou et~al.(2005)Chou, Chang, \protect\BIBand{}
  Lee}]{chou2005heuristic}
Chou FD, Chang TY, Lee CE (2005) A heuristic algorithm to minimize total
  weighted tardiness on a single machine with release times.
  \emph{International Transactions in Operational Research} 12(2):215--233.

\bibitem[{Della~Croce et~al.(1998)Della~Croce, Tadei, Baracco,
  \protect\BIBand{} Grosso}]{della1998new}
Della~Croce F, Tadei R, Baracco P, Grosso A (1998) A new decomposition approach
  for the single machine total tardiness scheduling problem. \emph{Journal of
  the Operational Research Society} 49(10):1101--1106.

\bibitem[{Ding et~al.(2016)Ding, L{\"u}, Cheng, \protect\BIBand{}
  Xu}]{ding2016breakout}
Ding J, L{\"u} Z, Cheng T, Xu L (2016) Breakout dynasearch for the
  single-machine total weighted tardiness problem. \emph{Computers \&
  Industrial Engineering} 98:1--10.

\bibitem[{Durasevi{\'c} \protect\BIBand{}
  Jakobovi{\'c}(2023)}]{ɖurasevic2023heuristic}
Durasevi{\'c} M, Jakobovi{\'c} D (2023) Heuristic and metaheuristic methods for
  the parallel unrelated machines scheduling problem: a survey.
  \emph{Artificial Intelligence Review} 56(4):3181--3289.

\bibitem[{Garey \protect\BIBand{} Johnson(1979)}]{garey-johnson-1979}
Garey MR, Johnson DS (1979) \emph{Computers and intractability}, volume 174
  (freeman San Francisco).

\bibitem[{Garraffa et~al.(2018)Garraffa, Shang, Della~Croce, \protect\BIBand{}
  T'kindt}]{garraffa2018exact}
Garraffa M, Shang L, Della~Croce F, T'kindt V (2018) An exact exponential
  branch-and-merge algorithm for the single machine total tardiness problem.
  \emph{Theoretical Computer Science} 745:133--149.

\bibitem[{Goemans et~al.(2002)Goemans, Queyranne, Schulz, Skutella,
  \protect\BIBand{} Wang}]{goemans2002single}
Goemans MX, Queyranne M, Schulz AS, Skutella M, Wang Y (2002) Single machine
  scheduling with release dates. \emph{SIAM Journal on Discrete Mathematics}
  15(2):165--192.

\bibitem[{Goffin \protect\BIBand{} Kiwiel(1999)}]{goffin1999convergence}
Goffin JL, Kiwiel KC (1999) Convergence of a simple subgradient level method.
  \emph{Mathematical Programming} 85(1):207--211.

\bibitem[{Goodfellow et~al.(2016)Goodfellow, Bengio, \protect\BIBand{}
  Courville}]{goodfellow2016deep}
Goodfellow I, Bengio Y, Courville A (2016) \emph{Deep learning} (MIT press).

\bibitem[{Graham et~al.(1979)Graham, Lawler, Lenstra, \protect\BIBand{}
  Kan}]{graham1979optimization}
Graham RL, Lawler EL, Lenstra JK, Kan AR (1979) Optimization and approximation
  in deterministic sequencing and scheduling: a survey. \emph{Annals of
  discrete mathematics}, volume~5, 287--326 (Elsevier).

\bibitem[{Grimes \protect\BIBand{} Hebrard(2015)}]{grimes2015solving}
Grimes D, Hebrard E (2015) Solving variants of the job shop scheduling problem
  through conflict-directed search. \emph{INFORMS Journal on Computing}
  27(2):268--284.

\bibitem[{Hearst et~al.(1998)Hearst, Dumais, Osuna, Platt, \protect\BIBand{}
  Scholkopf}]{hearst1998support}
Hearst MA, Dumais ST, Osuna E, Platt J, Scholkopf B (1998) Support vector
  machines. \emph{IEEE Intelligent Systems and their applications}
  13(4):18--28.

\bibitem[{Holsenback \protect\BIBand{} Russell(1992)}]{holsenback1992heuristic}
Holsenback JE, Russell RM (1992) A heuristic algorithm for sequencing on one
  machine to minimize total tardiness. \emph{Journal of the Operational
  Research Society} 43(1):53--62.

\bibitem[{Huang et~al.(2022)Huang, Wang, Zhou, \protect\BIBand{}
  Lu}]{huang2022planning}
Huang S, Wang Z, Zhou J, Lu J (2022) Planning irregular object packing via
  hierarchical reinforcement learning. \emph{IEEE Robotics and Automation
  Letters} 8(1):81--88.

\bibitem[{Janiak et~al.(2009)Janiak, Krysiak, Pappis, \protect\BIBand{}
  Voutsinas}]{janiak2009scheduling}
Janiak A, Krysiak T, Pappis CP, Voutsinas TG (2009) A scheduling problem with
  job values given as a power function of their completion times.
  \emph{European Journal of Operational Research} 193(3):836--848.

\bibitem[{Jouglet et~al.(2008)Jouglet, Savourey, Carlier, \protect\BIBand{}
  Baptiste}]{jouglet2008dominance}
Jouglet A, Savourey D, Carlier J, Baptiste P (2008) Dominance-based heuristics
  for one-machine total cost scheduling problems. \emph{European Journal of
  Operational Research} 184(3):879--899.

\bibitem[{Khalil et~al.(2016)Khalil, Le~Bodic, Song, Nemhauser,
  \protect\BIBand{} Dilkina}]{khalil2016learning}
Khalil E, Le~Bodic P, Song L, Nemhauser G, Dilkina B (2016) Learning to branch
  in mixed integer programming. \emph{Proceedings of the AAAI Conference on
  Artificial Intelligence}, volume~30.

\bibitem[{Kingma \protect\BIBand{} Ba(2014)}]{kingma2014adam}
Kingma DP, Ba J (2014) Adam: A method for stochastic optimization. \emph{arXiv
  preprint arXiv:1412.6980} .

\bibitem[{Kool et~al.(2019)Kool, van Hoof, \protect\BIBand{}
  Welling}]{kool2018attention}
Kool W, van Hoof H, Welling M (2019) Attention, learn to solve routing
  problems! \emph{International Conference on Learning Representations},
  \urlprefix\url{https://openreview.net/forum?id=ByxBFsRqYm}.

\bibitem[{LaValley(2008)}]{lavalley2008logistic}
LaValley MP (2008) Logistic regression. \emph{Circulation} 117(18):2395--2399.

\bibitem[{Lawler(1977{\natexlab{a}})}]{lawler-1977}
Lawler EL (1977{\natexlab{a}}) A “pseudopolynomial” algorithm for
  sequencing jobs to minimize total tardiness. \emph{Annals of discrete
  Mathematics}, volume~1, 331--342 (Elsevier).

\bibitem[{Lawler(1977{\natexlab{b}})}]{lawler1977pseudopolynomial}
Lawler EL (1977{\natexlab{b}}) A “pseudopolynomial” algorithm for
  sequencing jobs to minimize total tardiness. \emph{Annals of discrete
  Mathematics}, volume~1, 331--342 (Elsevier).

\bibitem[{LeCun et~al.(2015)LeCun, Bengio, \protect\BIBand{}
  Hinton}]{lecun2015deep}
LeCun Y, Bengio Y, Hinton G (2015) Deep learning. \emph{nature}
  521(7553):436--444.

\bibitem[{Lecun et~al.(1998)Lecun, Bottou, Bengio, \protect\BIBand{}
  Haffner}]{lenet1998}
Lecun Y, Bottou L, Bengio Y, Haffner P (1998) Gradient-based learning applied
  to document recognition. \emph{Proceedings of the IEEE} 86(11):2278--2324,
  \urlprefix\url{http://dx.doi.org/10.1109/5.726791}.

\bibitem[{Lenstra et~al.(1977)Lenstra, Kan, \protect\BIBand{}
  Brucker}]{lenstra1977complexity}
Lenstra JK, Kan AR, Brucker P (1977) Complexity of machine scheduling problems.
  \emph{Annals of discrete mathematics}, volume~1, 343--362 (Elsevier).

\bibitem[{Li et~al.(2023)Li, Lang, Hong, \protect\BIBand{}
  Reggelin}]{li2023two}
Li F, Lang S, Hong B, Reggelin T (2023) A two-stage rnn-based deep
  reinforcement learning approach for solving the parallel machine scheduling
  problem with due dates and family setups. \emph{Journal of Intelligent
  Manufacturing} 1--34.

\bibitem[{Li et~al.(2022)Li, Liu, Yang, Peng, \protect\BIBand{}
  Zhou}]{CNNserveyZeWenLi}
Li Z, Liu F, Yang W, Peng S, Zhou J (2022) A survey of convolutional neural
  networks: Analysis, applications, and prospects. \emph{IEEE Transactions on
  Neural Networks and Learning Systems} 33(12):6999--7019,
  \urlprefix\url{http://dx.doi.org/10.1109/TNNLS.2021.3084827}.

\bibitem[{Liu et~al.(2023)Liu, Luh, Sun, Bragin, \protect\BIBand{}
  Yan}]{liu2023integrating}
Liu A, Luh PB, Sun K, Bragin MA, Yan B (2023) Integrating machine learning and
  mathematical optimization for job shop scheduling. \emph{IEEE Transactions on
  Automation Science and Engineering} .

\bibitem[{Minaeva et~al.(2016)Minaeva, {\v{S}}{\u{u}}cha, Akesson,
  \protect\BIBand{} Hanz{\'a}lek}]{minaeva2016scalable}
Minaeva A, {\v{S}}{\u{u}}cha P, Akesson B, Hanz{\'a}lek Z (2016) Scalable and
  efficient configuration of time-division multiplexed resources. \emph{Journal
  of Systems and Software} 113:44--58.

\bibitem[{Minella et~al.(2008)Minella, Ruiz, \protect\BIBand{}
  Ciavotta}]{minella2008review}
Minella G, Ruiz R, Ciavotta M (2008) A review and evaluation of multiobjective
  algorithms for the flowshop scheduling problem. \emph{INFORMS Journal on
  Computing} 20(3):451--471.

\bibitem[{Moore(1968)}]{moore-1968}
Moore JM (1968) An n job, one machine sequencing algorithm for minimizing the
  number of late jobs. \emph{Management science} 15(1):102--109.

\bibitem[{Morabit et~al.(2021)Morabit, Desaulniers, \protect\BIBand{}
  Lodi}]{morabit2021machine}
Morabit M, Desaulniers G, Lodi A (2021) Machine-learning--based column
  selection for column generation. \emph{Transportation Science}
  55(4):815--831.

\bibitem[{Morabit et~al.(2023)Morabit, Desaulniers, \protect\BIBand{}
  Lodi}]{morabit2023machine}
Morabit M, Desaulniers G, Lodi A (2023) Machine-learning--based arc selection
  for constrained shortest path problems in column generation. \emph{INFORMS
  Journal on Optimization} 5(2):191--210.

\bibitem[{Murtagh(1991)}]{murtagh1991multilayer}
Murtagh F (1991) Multilayer perceptrons for classification and regression.
  \emph{Neurocomputing} 2(5-6):183--197.

\bibitem[{Nazari et~al.(2018)Nazari, Oroojlooy, Snyder, \protect\BIBand{}
  Tak{\'a}c}]{nazari2018reinforcement}
Nazari M, Oroojlooy A, Snyder L, Tak{\'a}c M (2018) Reinforcement learning for
  solving the vehicle routing problem. \emph{Advances in neural information
  processing systems} 31.

\bibitem[{Nedic \protect\BIBand{} Bertsekas(2001)}]{nedic2001incremental}
Nedic A, Bertsekas DP (2001) Incremental subgradient methods for
  nondifferentiable optimization. \emph{SIAM Journal on Optimization}
  12(1):109--138.

\bibitem[{Parmentier \protect\BIBand{}
  T’kindt(2023)}]{parmentier2023structured}
Parmentier A, T’kindt V (2023) Structured learning based heuristics to solve
  the single machine scheduling problem with release times and sum of
  completion times. \emph{European Journal of Operational Research}
  305(3):1032--1041.

\bibitem[{Paulus \protect\BIBand{} Krause(2024)}]{paulus2024learning}
Paulus M, Krause A (2024) Learning to dive in branch and bound. \emph{Advances
  in Neural Information Processing Systems} 36.

\bibitem[{Pinedo(2012)}]{pinedo2012scheduling}
Pinedo ML (2012) \emph{Scheduling}, volume~29 (Springer).

\bibitem[{Polyak(1969)}]{polyak1969minimization}
Polyak BT (1969) Minimization of unsmooth functionals. \emph{USSR Computational
  Mathematics and Mathematical Physics} 9(3):14--29.

\bibitem[{Potts \protect\BIBand{} Van~Wassenhove(1991)}]{potts1991single}
Potts C, Van~Wassenhove LN (1991) Single machine tardiness sequencing
  heuristics. \emph{IIE transactions} 23(4):346--354.

\bibitem[{Rumelhart et~al.(1986)Rumelhart, Hinton, \protect\BIBand{}
  Williams}]{rumelhart1986learning}
Rumelhart DE, Hinton GE, Williams RJ (1986) Learning representations by
  back-propagating errors. \emph{nature} 323(6088):533--536.

\bibitem[{Sahni(1976)}]{sahni-1976}
Sahni SK (1976) Algorithms for scheduling independent tasks. \emph{Journal of
  the ACM (JACM)} 23(1):116--127.

\bibitem[{Schmidt \protect\BIBand{} Stober(2021)}]{schmidt2021approaching}
Schmidt J, Stober S (2021) Approaching scheduling problems via a deep hybrid
  greedy model and supervised learning. \emph{IFAC-PapersOnLine}
  54(1):805--810.

\bibitem[{Shen et~al.(2022)Shen, Sun, Li, Eberhard, \protect\BIBand{}
  Ernst}]{shen2022enhancing}
Shen Y, Sun Y, Li X, Eberhard A, Ernst A (2022) Enhancing column generation by
  a machine-learning-based pricing heuristic for graph coloring.
  \emph{Proceedings of the AAAI conference on artificial intelligence},
  volume~36, 9926--9934.

\bibitem[{Smith et~al.(1956)}]{smith-1956}
Smith WE, et~al. (1956) Various optimizers for single-stage production.
  \emph{Naval Research Logistics Quarterly} 3(1-2):59--66.

\bibitem[{Sousa \protect\BIBand{} Wolsey(1992)}]{sousa1992time}
Sousa JP, Wolsey LA (1992) A time indexed formulation of non-preemptive single
  machine scheduling problems. \emph{Mathematical Programming} 54:353--367.

\bibitem[{Szwarc et~al.(1988)Szwarc, Posner, \protect\BIBand{}
  Liu}]{szwarc1988single}
Szwarc W, Posner ME, Liu JJ (1988) The single machine problem with a quadratic
  cost function of completion times. \emph{Management Science}
  34(12):1480--1488.

\bibitem[{Tasgetiren et~al.(2004)Tasgetiren, Sevkli, Liang, \protect\BIBand{}
  Gen{\c{c}}yilmaz}]{tasgetiren2004particle}
Tasgetiren MF, Sevkli M, Liang YC, Gen{\c{c}}yilmaz G (2004) Particle swarm
  optimization algorithm for single machine total weighted tardiness problem.
  \emph{Proceedings of the 2004 Congress on Evolutionary Computation},
  volume~2, 1412--1419 (IEEE).

\bibitem[{T'kindt \protect\BIBand{} Billaut(2001)}]{t2001multicriteria}
T'kindt V, Billaut JC (2001) Multicriteria scheduling problems: a survey.
  \emph{RAIRO-Operations Research} 35(2):143--163.

\bibitem[{Uchro{\'n}ski(2021)}]{uchronski2021parallel}
Uchro{\'n}ski M (2021) Parallel algorithm with blocks for a single-machine
  total weighted tardiness scheduling problem. \emph{Applied Sciences}
  11(5):2069.

\bibitem[{V{\'a}clav{\'\i}k et~al.(2018)V{\'a}clav{\'\i}k, Nov{\'a}k,
  {\v{S}}{\u{u}}cha, \protect\BIBand{} Hanz{\'a}lek}]{vaclavik2018accelerating}
V{\'a}clav{\'\i}k R, Nov{\'a}k A, {\v{S}}{\u{u}}cha P, Hanz{\'a}lek Z (2018)
  Accelerating the branch-and-price algorithm using machine learning.
  \emph{European Journal of Operational Research} 271(3):1055--1069.

\bibitem[{Van~den Akker et~al.(2000)Van~den Akker, Hurkens, \protect\BIBand{}
  Savelsbergh}]{van2000time}
Van~den Akker J, Hurkens CA, Savelsbergh MW (2000) Time-indexed formulations
  for machine scheduling problems: Column generation. \emph{INFORMS Journal on
  Computing} 12(2):111--124.

\bibitem[{Vaswani et~al.(2017)Vaswani, Shazeer, Parmar, Uszkoreit, Jones,
  Gomez, Kaiser, \protect\BIBand{} Polosukhin}]{vaswani2017attention}
Vaswani A, Shazeer N, Parmar N, Uszkoreit J, Jones L, Gomez AN, Kaiser {\L},
  Polosukhin I (2017) Attention is all you need. \emph{Advances in neural
  information processing systems} 30.

\bibitem[{Vinyals et~al.(2015)Vinyals, Fortunato, \protect\BIBand{}
  Jaitly}]{vinyals2015pointer}
Vinyals O, Fortunato M, Jaitly N (2015) Pointer networks. \emph{Advances in
  neural information processing systems} 28.

\bibitem[{Wang \protect\BIBand{} Usher(2005)}]{wang2005application}
Wang YC, Usher JM (2005) Application of reinforcement learning for agent-based
  production scheduling. \emph{Engineering applications of artificial
  intelligence} 18(1):73--82.

\bibitem[{Weckman et~al.(2008)Weckman, Ganduri, \protect\BIBand{}
  Koonce}]{weckman2008neural}
Weckman GR, Ganduri CV, Koonce DA (2008) A neural network job-shop scheduler.
  \emph{Journal of Intelligent Manufacturing} 19:191--201.

\bibitem[{Yuan et~al.(2016)Yuan, Jiang, \protect\BIBand{}
  Wang}]{yuan2016dynamic}
Yuan B, Jiang Z, Wang L (2016) Dynamic parallel machine scheduling with random
  breakdowns using the learning agent. \emph{International Journal of Services
  Operations and Informatics} 8(2):94--103.

\bibitem[{Zhang et~al.(2020)Zhang, Song, Cao, Zhang, Tan, \protect\BIBand{}
  Chi}]{zhang2020learning}
Zhang C, Song W, Cao Z, Zhang J, Tan PS, Chi X (2020) Learning to dispatch for
  job shop scheduling via deep reinforcement learning. \emph{Advances in Neural
  Information Processing Systems} 33:1621--1632.

\bibitem[{Zhang et~al.(2022)Zhang, Chen, Gendreau, \protect\BIBand{}
  Langevin}]{zhang2022learning}
Zhang X, Chen L, Gendreau M, Langevin A (2022) Learning-based branch-and-price
  algorithms for the vehicle routing problem with time windows and
  two-dimensional loading constraints. \emph{INFORMS Journal on Computing}
  34(3):1419--1436.

\end{thebibliography}

\clearpage

\section*{Appendix A. Background Knowledge of the Backpropagation Algorithm}
\label{sec:back_propagation}
\addcontentsline{toc}{section}{Appendix A}

The backpropagration algorithm \citep{rumelhart1986learning} is used to calculate the gradient for offline supervised training and online single-instance learning. To help readers understand our paper, in this appendix, we use the computation of the gradient of a multilayer perceptron (MLP) as an example to briefly introduce the main concepts of the backpropagation algorithm.
An MLP consists of fully connected neurons and nonlinear activation functions \citep{murtagh1991multilayer}. For the $L$-layer MLP, given a vector of input $x=x^{(1)}\in \mathbb{R}^{d_1}$, affine mappings and non-linear activation functions are iteratively applied to generate the subsequent outputs, i.e., $x^{(l+1)}=\sigma(W^{(l)} x^{(l)} + b^{(l)})$, for $l=1,\cdots,L$, where $W^{(l)}\in\mathbb{R}^{d_l\times d_{l+1}}$ and $b^{(l)}\in\mathbb{R}^{d_{l+1}}$ are the learnable parameters (called weights and biases) of the $l$-th layer, $\sigma(\cdot)$ is an activation function, and $x^{(L+1)}$ is the final output of the MLP. 

When an MLP has a large number of layers, there will be a large number of learnable parameters and many non-linear activation functions, and the gradient-based algorithm is almost the only reliable algorithm for training. Therefore, the key to training lies in computing the gradient of loss value $\mathcal{L}(y,f_\theta(x))$ with respect to the parameters $\theta$. 
By leveraging the chain rule of differentiation, \citet{rumelhart1986learning} developed the famous backpropagation algorithm, where the gradients are calculated layer by layer, iterating backward from the last layer. For example, for the $L$-layer MLP described above, the gradients of the loss with respect to the weights and biases of the $l$th layer can be calculated as follows:
\begin{equation*}
\begin{aligned}
    \nabla_{W^{(l)}} \mathcal{L}(y,f_\theta(x)) = D^{\top}_{W^{(l)}} x^{(l+1)} \nabla_{x^{(l+1)}}  \mathcal{L}(y,f_\theta(x)),\\
    \nabla_{b^{(l)}} \mathcal{L}(y,f_\theta(x)) = D^{\top}_{b^{(l)}} x^{(l+1)} \nabla_{x^{(l+1)}}  \mathcal{L}(y,f_\theta(x)),
\end{aligned}
\end{equation*}
with $\nabla_{x^{(l+1)}} \mathcal{L}(y,f_\theta(x)) = D^{\top}_{x^{(l+1)}} x^{(l+2)} \nabla_{x^{(l+2)}}  \mathcal{L}(y,f_\theta(x))$, where $D$ denotes the Jacobian matrix with respect to the subscript. The Jacobian operator recurring throughout the recursion indicates that the entire model/computation process has to be differentiable from the input $x^{(1)}$ to the output $x^{(L+1)}$.
\begin{figure}[h]
    \centering
    \includegraphics[width=0.4\linewidth]{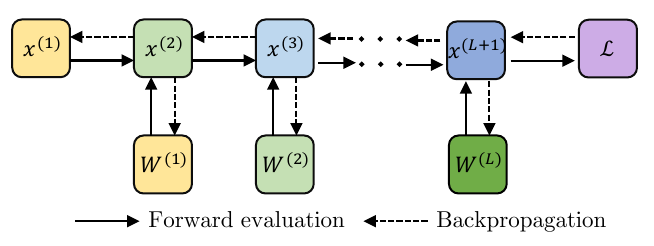}
    \caption{Computation flows in MLPs}
    \label{fig:bbvlm}
\end{figure}

\section*{Appendix B. Supplementary Materials for Online Single-instance Learning }
\label{sec:back_propagation}
\addcontentsline{toc}{section}{Appendix B}

\subsection{Pseudocode}
We summarize the steps of the online single instance learning as Algorithm \ref{alg:online_learning}. 
\begin{algorithm}
    \caption{Online Single-instance Learning}
    \label{alg:online_learning}
    \begin{algorithmic}[1]
    \STATE {\bfseries Input:} Hyper-parameters $\varphi_0$ and $\Delta\varphi$ for the smooth feasibility surrogate, learnable parameters $\theta^{\mathrm{offline}}$, hyper-parameter $\varepsilon$ and $B$ for adjusting the target, and hyper-parameter $B^{\mathrm{stop}}$ for stopping criteria. 
       \STATE {\bfseries Initialize:} $\theta_0\leftarrow \theta^{\mathrm{offline}}$, $\lambda_0\leftarrow \varepsilon\cdot J(\theta^{\mathrm{offline}})$,
       $\theta^{\mathrm{best}}\leftarrow \theta^{\mathrm{offline}}$.
       \REPEAT
          \IF{$J(\theta_k)<J(\theta^{\mathrm{best}})$}
            \STATE $\theta^{\mathrm{best}} \leftarrow \theta_k$.
          \ENDIF
          \STATE  $g_k \leftarrow \nabla_{\theta} \tilde{J}(\theta_k)$.
          \STATE $\Bar{J}_{k} \leftarrow \inf_{\kappa<k} \tilde{J}(\theta_\kappa) - \lambda_k$.
          \STATE $s_k\leftarrow \frac{\tilde{J}(\theta_k)-\Bar{J}_{k}}{||\nabla_{\theta} \tilde{J}(\theta_k)||_2^2}$.
          \STATE $\theta_{k+1} \leftarrow \theta_k-s_k g_k$.
          \IF{$k-\argmin_{\kappa} \tilde{J}(\theta_{\kappa})> B$}
            \STATE $\lambda_{k+1} \leftarrow \lambda_{k} /2 $, $\varphi_{k+1} \leftarrow \varphi_{k} + \Delta\varphi$,
            $\theta_{k+1}\leftarrow \theta^{\mathrm{best}}  $.
            \ELSE
            \STATE $\lambda_{k+1} \leftarrow \lambda_{k}$, $\varphi_{k+1} \leftarrow \varphi_{k}$,
          \ENDIF
        \UNTIL  stopping criteria $k-\argmin_{\kappa} J(\theta_\kappa)>B^{\mathrm{stop}}$ is satisfied. 
    \STATE $\theta^{\mathrm{online}}\leftarrow \theta^{k}$.
    \STATE {\bfseries Output:} Trained parameters $\theta^{\mathrm{online}} $, job starting times $ F_{\mathcal{V}}(G_{\mathcal{I}}(\theta^{\mathrm{online}}))$
    \end{algorithmic}
\end{algorithm}

\subsection{Hyperparameters}
With the learnable parameters obtained offline, 
the online single-instance learning approach is used to find a solution with an improved quality for each online testing instance. The hyperparameters used are exactly the same for all the problem types. We set $B=50$,  $B^{\mathrm{stop}}=100$,  $\varepsilon=0.02$, $\varphi_0=0.5$ and $\Delta\varphi=0.1$.

\subsection{Illustration of the Online Single-instance Learning}
\label{illustration_online}

To visualize the online single-instance learning, we show the values of $\{J(\theta_k)\}$ and $\{\tilde{J}(\theta_k)\}$ obtained when solving an instance of $1|r_j|\sum w_j T_j$ with 700 jobs. Since the step functions are approximated by sigmoid functions within the feasibility surrogate, there is an error gap between $\tilde{J(\theta_k)}$ and the true objective value $J(\theta_k)$. However, as $\tilde{J(\theta_k)}$ being optimized, the true objective value $J(\theta_k)$ is also optimized.

\begin{figure}[!t]
\begin{center}
        \includegraphics[width=0.7\linewidth]{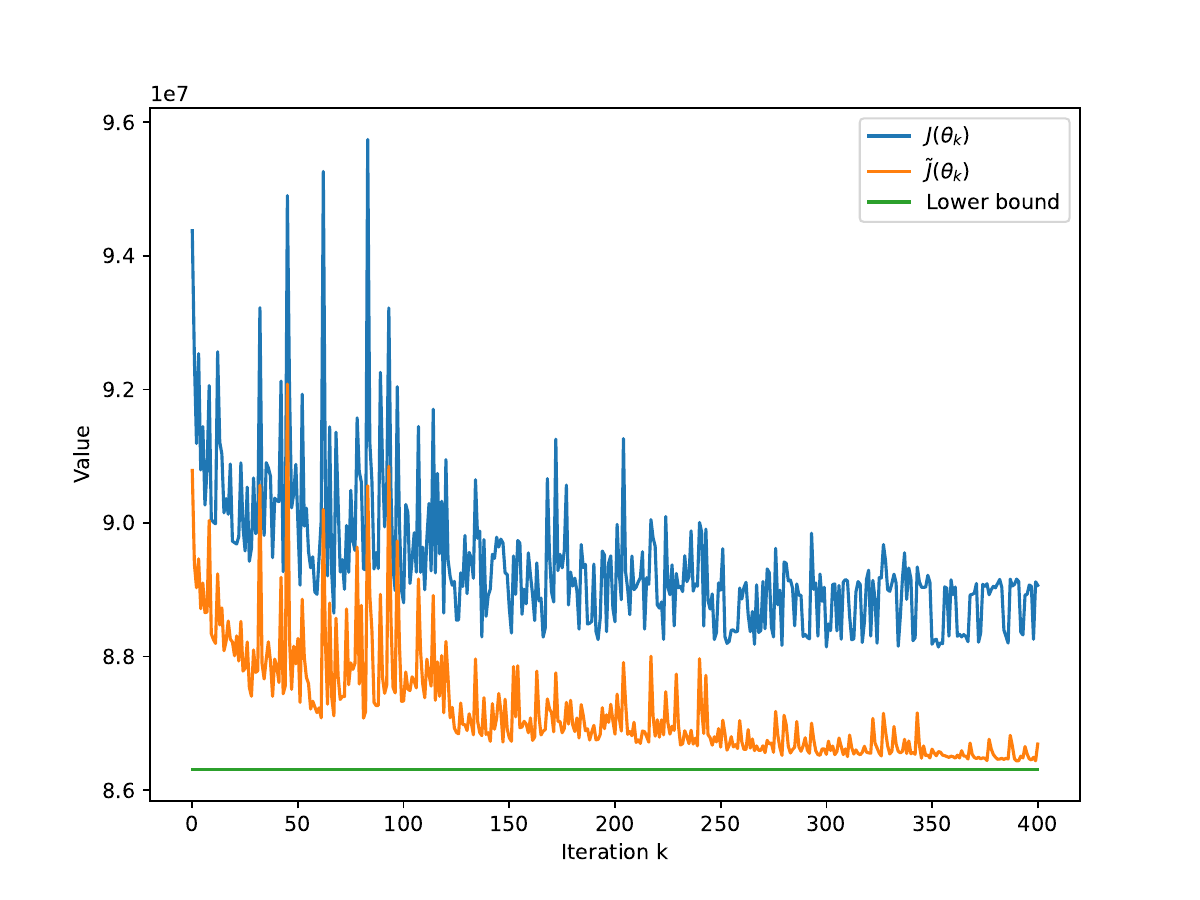}
    \end{center}
    \caption{Online single-instance learning for an instance of $1|r_j|\sum w_jT_j$ with 700 jobs}
    \label{Fig:online_loss}
\end{figure}

\section*{Appendix C. Layer sizes of the UMSNN used in the numerical result section}
\label{layer_sizes}
\addcontentsline{toc}{section}{Appendix C}
In this section, we present the layer sizes of the UMSNN used in the numerical testing section. In the numerical results, two problem sizes are considered and are solved with separate UMSNN models. 
The layer sizes for the two different problem sizes are slightly different, and are presented separately in the following.

\subsection{Layer sizes for instances with 500 to 700 jobs}
\textbf{Sizes of the input and output:}
Each job’s processing time is at most 100, so in  \eqref{eq:coding_pros}, $\Tilde{p_j}$ is encoded as a length-100 vector. The total number of time slots (i.e., $T$) to 53000,  which is sufficiently large for all instances considered. Hence, each job has 53,000 starting costs. The release date encoding $\tilde{r}_j$ also has a length of 53,000. We group every 100 time slots into one time window, yielding 340 time windows in total. Therefore, in \eqref{UMSNN_output}, $\gamma=340$.

\textbf{Size of the input module:}
The MLP in the input module for processing $\Tilde{p_j}$ is formed by stacking a linear layer with a size of $100 \times 256$, a ReLU activation function, and a linear layer with a size of $256 \times 256$. For ease of representation, we denote this architecture as [Linear(100, 256), ReLU, Linear(256, 256)]. For Example 1, the CNN in the input module for processing the starting costs has a structure of [Conv1d(9, 32, 3), ReLU, Conv1d(6, 64, 2), ReLU, Conv1d(3, 64, 2), ReLU, Conv1d(1, 1, 1), ReLU, Linear($2831\times1024$)], and the CNN for release dates has a structure of [Conv1d(9, 8, 3), ReLU, Conv1d(6, 16, 2), ReLU, Conv1d(3, 8, 2), ReLU, Conv1d(1, 1, 1), ReLU, Linear($2831\times1024$)]. 
Here, Conv1d(9, 32, 3) denotes a 1-D CNN layer, which has 32 kernels each of size 9, and a stride of 3.

\textbf{Size of the encoder:}
We employ a transformer-based encoder with 8 layers, each having 8 heads. The feedforward layer in each transformer block has dimension 2048.

\textbf{Size of the output module:}
The output layer has a structure of [Linear($1280 \times 1280$), ReLU, Linear($1280 \times 340$)].

\subsection{Layer sizes for instances with 800 to 1000 jobs}

When learning to solve instances with 800 to 1000 jobs, the total number of time slots (i.e., $T$) is set as 68000. Also, we group every 100 time slots into a single window, resulting in 680 time windows. The layers in the input module are identical to those for the 500–700-job case, except that the last linear layer is 
Linear($5664\times1024$). The encoder remains the same as in the 500–700-job setting. The output layer has a structure of [Linear($1280 \times 1280$), ReLU, Linear($1280 \times 680$)].

\end{document}